\documentstyle[amssymb]{article}

\newenvironment{align}{\[}{\]}
\def\address#1{address: #1\par}
\def\email#1{email: #1\par}
\def\subjclass#1{Ams classification: #1\par}

\newcommand{\rest}{\restriction}

\makeatletter
\let\DOTSI\relax
\def\RIfM@{\relax\ifmmode}%
\def\FN@{\futurelet\next}%
\newcount\intno@
\def\iint{\DOTSI\intno@\tw@\FN@\ints@}%
\def\iiint{\DOTSI\intno@\thr@@\FN@\ints@}%
\def\iiiint{\DOTSI\intno@4 \FN@\ints@}%
\def\idotsint{\DOTSI\intno@\z@\FN@\ints@}%
\def\ints@{\findlimits@\ints@@}%
\newif\iflimtoken@
\newif\iflimits@
\def\findlimits@{\limtoken@true\ifx\next\limits\limits@true
 \else\ifx\next\nolimits\limits@false\else
 \limtoken@false\ifx\ilimits@\nolimits\limits@false\else
 \ifinner\limits@false\else\limits@true\fi\fi\fi\fi}%
\def\multint@{\int\ifnum\intno@=\z@\intdots@                                
 \else\intkern@\fi                                                          
 \ifnum\intno@>\tw@\int\intkern@\fi                                         
 \ifnum\intno@>\thr@@\int\intkern@\fi                                       
 \int}
\def\multintlimits@{\intop\ifnum\intno@=\z@\intdots@\else\intkern@\fi
 \ifnum\intno@>\tw@\intop\intkern@\fi
 \ifnum\intno@>\thr@@\intop\intkern@\fi\intop}%
\def\intic@{\mathchoice{\hskip.5em}{\hskip.4em}{\hskip.4em}{\hskip.4em}}%
\def\negintic@{\mathchoice
 {\hskip-.5em}{\hskip-.4em}{\hskip-.4em}{\hskip-.4em}}%
\def\ints@@{\iflimtoken@                                                    
 \def\ints@@@{\iflimits@\negintic@\mathop{\intic@\multintlimits@}\limits    
  \else\multint@\nolimits\fi                                                
  \eat@}
 \else                                                                      
 \def\ints@@@{\iflimits@\negintic@
  \mathop{\intic@\multintlimits@}\limits\else
  \multint@\nolimits\fi}\fi\ints@@@}%
\def\intkern@{\mathchoice{\!\!\!}{\!\!}{\!\!}{\!\!}}%
\def\plaincdots@{\mathinner{\cdotp\cdotp\cdotp}}%
\def\intdots@{\mathchoice{\plaincdots@}%
 {{\cdotp}\mkern1.5mu{\cdotp}\mkern1.5mu{\cdotp}}%
 {{\cdotp}\mkern1mu{\cdotp}\mkern1mu{\cdotp}}%
 {{\cdotp}\mkern1mu{\cdotp}\mkern1mu{\cdotp}}}%
\def\rmfam{\z@}%
\newif\iffirstchoice@
\firstchoice@true
\def\textfonti{\the\textfont\@ne}%
\def\textfontii{\the\textfont\tw@}%
\def\text{\RIfM@\expandafter\text@\else\expandafter\text@@\fi}%
\def\text@@#1{\leavevmode\hbox{#1}}%
\def\text@#1{\mathchoice
 {\hbox{\everymath{\displaystyle}\def\textfonti{\the\textfont\@ne}%
  \def\textfontii{\the\textfont\tw@}\textdef@@ T#1}}%
 {\hbox{\firstchoice@false
  \everymath{\textstyle}\def\textfonti{\the\textfont\@ne}%
  \def\textfontii{\the\textfont\tw@}\textdef@@ T#1}}%
 {\hbox{\firstchoice@false
  \everymath{\scriptstyle}\def\textfonti{\the\scriptfont\@ne}%
  \def\textfontii{\the\scriptfont\tw@}\textdef@@ S\rm#1}}%
 {\hbox{\firstchoice@false
  \everymath{\scriptscriptstyle}\def\textfonti
  {\the\scriptscriptfont\@ne}%
  \def\textfontii{\the\scriptscriptfont\tw@}\textdef@@ s\rm#1}}}%
\def\textdef@@#1{\textdef@#1\rm\textdef@#1\bf\textdef@#1\sl\textdef@#1\it}%
\def\DN@{\def\next@}%
\def\eat@#1{}%
\def\textdef@#1#2{%
 \DN@{\csname\expandafter\eat@\string#2fam\endcsname}%
 \if S#1\edef#2{\the\scriptfont\next@\relax}%
 \else\if s#1\edef#2{\the\scriptscriptfont\next@\relax}%
 \else\edef#2{\the\textfont\next@\relax}\fi\fi}%
\def\Let@{\relax\iffalse{\fi\let\\=\cr\iffalse}\fi}%
\def\vspace@{\def\vspace##1{\crcr\noalign{\vskip##1\relax}}}%
\def\multilimits@{\bgroup\vspace@\Let@
 \baselineskip\fontdimen10 \scriptfont\tw@
 \advance\baselineskip\fontdimen12 \scriptfont\tw@
 \lineskip\thr@@\fontdimen8 \scriptfont\thr@@
 \lineskiplimit\lineskip
 \vbox\bgroup\ialign\bgroup\hfil$\m@th\scriptstyle{##}$\hfil\crcr}%
\def\Sb{_\multilimits@}%
\def\endSb{\crcr\egroup\egroup\egroup}%
\def\Sp{^\multilimits@}%

\newdimen\ex@
\def\rightarrowfill@#1{$#1\m@th\mathord-\mkern-6mu\cleaders
 \hbox{$#1\mkern-2mu\mathord-\mkern-2mu$}\hfill
 \mkern-6mu\mathord\rightarrow$}%
\def\leftarrowfill@#1{$#1\m@th\mathord\leftarrow\mkern-6mu\cleaders
 \hbox{$#1\mkern-2mu\mathord-\mkern-2mu$}\hfill\mkern-6mu\mathord-$}%
\def\leftrightarrowfill@#1{$#1\m@th\mathord\leftarrow\mkern-6mu\cleaders
 \hbox{$#1\mkern-2mu\mathord-\mkern-2mu$}\hfill
 \mkern-6mu\mathord\rightarrow$}%
\def\overrightarrow{\mathpalette\overrightarrow@}%
\def\overrightarrow@#1#2{\vbox{\ialign{##\crcr\rightarrowfill@#1\crcr
 \noalign{\kern-\ex@\nointerlineskip}$\m@th\hfil#1#2\hfil$\crcr}}}%

\def\overleftarrow{\mathpalette\overleftarrow@}%
\def\overleftarrow@#1#2{\vbox{\ialign{##\crcr\leftarrowfill@#1\crcr
 \noalign{\kern-\ex@\nointerlineskip}$\m@th\hfil#1#2\hfil$\crcr}}}%
\def\overleftrightarrow{\mathpalette\overleftrightarrow@}%
\def\overleftrightarrow@#1#2{\vbox{\ialign{##\crcr\leftrightarrowfill@#1\crcr
 \noalign{\kern-\ex@\nointerlineskip}$\m@th\hfil#1#2\hfil$\crcr}}}%
\def\underrightarrow{\mathpalette\underrightarrow@}%
\def\underrightarrow@#1#2{\vtop{\ialign{##\crcr$\m@th\hfil#1#2\hfil$\crcr
 \noalign{\nointerlineskip}\rightarrowfill@#1\crcr}}}%

\def\underleftarrow{\mathpalette\underleftarrow@}%
\def\underleftarrow@#1#2{\vtop{\ialign{##\crcr$\m@th\hfil#1#2\hfil$\crcr
 \noalign{\nointerlineskip}\leftarrowfill@#1\crcr}}}%
\def\underleftrightarrow{\mathpalette\underleftrightarrow@}%
\def\underleftrightarrow@#1#2{\vtop{\ialign{##\crcr$\m@th\hfil#1#2\hfil$\crcr
 \noalign{\nointerlineskip}\leftrightarrowfill@#1\crcr}}}%
\newcount\GRAPHICSTYPE
\GRAPHICSTYPE=\z@
\def\GRAPHICSPS#1{%
 \ifnum\GRAPHICSTYPE=\@ne language "PS", include "#1"\else ps: #1\fi
 }%
\def\graffile#1#2#3#4{%
 \leavevmode\raise -#4 \hbox{%
  \raise #3 \hbox{\rule{0.003in}{0.003in}\special{#1}}%
  }%
 {\raise -#4 \hbox to #2 {\vrule height#3 width\z@ depth\z@\hfil}}%
 }%
\def\draftbox#1#2#3#4{%
 \leavevmode\raise -#4 \hbox{%
  \frame{\rlap{\protect\tiny #1}\hbox to #2%
   {\vrule height#3 width\z@ depth\z@\hfil}%
  }%
 }%
}%
\newcount\draft
\draft=\z@
\def\GRAPHIC#1#2#3#4#5{%
 \ifnum\draft=\@ne \draftbox{#2}{#3}{#4}{#5}%
  \else \graffile{#1}{#3}{#4}{#5}%
  \fi
 }%
\def\addtoLaTeXparams#1{\edef\LaTeXparams{\LaTeXparams #1}}%
\def\doFRAMEparams#1{\readFRAMEparams#1\end}%
\def\readFRAMEparams#1{%
 \ifx#1\end%
  \let\next=\relax
  \else
  \ifx#1i\dispkind=\z@\fi
  \ifx#1d\dispkind=\@ne\fi
  \ifx#1f\dispkind=\tw@\fi
  \ifx#1t\addtoLaTeXparams{t}\fi
  \ifx#1b\addtoLaTeXparams{b}\fi
  \ifx#1p\addtoLaTeXparams{p}\fi
  \ifx#1h\addtoLaTeXparams{h}\fi
  \let\next=\readFRAMEparams
  \fi
 \next
 }%
\def\IFRAME#1#2#3#4#5{\GRAPHIC{#5}{#4}{#1}{#2}{#3}}%
\def\DFRAME#1#2#3#4{%
 \begin{center}\GRAPHIC{#4}{#3}{#1}{#2}{\z@}\end{center}%
 }%
\def\FFRAME#1#2#3#4#5#6#7{%
 \begin{figure}[#1]%
  \begin{center}\GRAPHIC{#7}{#6}{#2}{#3}{\z@}\end{center}%
  \caption{\label{#5}#4}%
  \end{figure}%
 }%
\newcount\dispkind%
\def\FRAME#1#2#3#4#5#6#7#8{%
 \def\LaTeXparams{}%
 \dispkind=\z@
 \def\LaTeXparams{}%
 \doFRAMEparams{#1}%
 \ifnum\dispkind=\z@\IFRAME{#2}{#3}{#4}{#7}{#8}\else
  \ifnum\dispkind=\@ne\DFRAME{#2}{#3}{#7}{#8}\else
   \ifnum\dispkind=\tw@
    \edef\@tempa{\noexpand\FFRAME{\LaTeXparams}}%
    \@tempa{#2}{#3}{#5}{#6}{#7}{#8}%
    \fi
   \fi
  \fi
 }%
\def\limfunc#1{\mbox{\rm #1}}%
\long\def\QQQ#1#2{\long\expandafter\def\csname#1\endcsname{#2}}%
\def\QTP#1{}%
\long\def\QQA#1#2{}%
\def\QTR#1#2{{\csname#1\endcsname #2}}
\def\EXPAND#1[#2]#3{}%
\def\NOEXPAND#1[#2]#3{}%
\def\LaTeXparent#1{}%
\def\QTagDef#1#2#3{}%
\def\QQfnmark#1{\footnotemark}

\@ifundefined{INDEX}{\def\INDEX#1#2{}{}}{}%
\@ifundefined{SUBINDEX}{\def\SUBINDEX#1#2#3{}{}{}}{}%
\def\initial#1{\bigbreak{\raggedright\large\bf #1}\kern 2\p@\penalty3000}%
\@ifundefined{abstract}{%
 \def\abstract{%
  \if@twocolumn
   \section*{Abstract (Not appropriate in this style!)}%
   \else \small 
   \begin{center}{\bf Abstract\vspace{-.5em}\vspace{\z@}}\end{center}%
   \quotation 
   \fi
  }%
 }{%
 }%
\@ifundefined{endabstract}{\def\endabstract
  {\if@twocolumn\else\endquotation\fi}}{}%
\@ifundefined{maketitle}{\def\maketitle#1{}}{}%
\@ifundefined{affiliation}{\def\affiliation#1{}}{}%
\@ifundefined{proof}{}{}%
\@ifundefined{endproof}{}{}%
\@ifundefined{newfield}{\def\newfield#1#2{}}{}%
\@ifundefined{chapter}{\def\chapter#1{\par(Chapter head:)#1\par }%
 \newcount\c@chapter}{}%
\@ifundefined{part}{\def\part#1{\par(Part head:)#1\par }}{}%
\@ifundefined{section}{\def\section#1{\par(Section head:)#1\par }}{}%
\@ifundefined{subsection}{\def\subsection#1%
 {\par(Subsection head:)#1\par }}{}%
\@ifundefined{subsubsection}{\def\subsubsection#1%
 {\par(Subsubsection head:)#1\par }}{}%
\@ifundefined{paragraph}{\def\paragraph#1%
 {\par(Subsubsubsection head:)#1\par }}{}%
\@ifundefined{subparagraph}{\def\subparagraph#1%
 {\par(Subsubsubsubsection head:)#1\par }}{}%
\@ifundefined{therefore}{}{}%
\@ifundefined{backepsilon}{}{}%
\@ifundefined{yen}{}{}%
\@ifundefined{registered}{\def\registered{\relax\ifmmode{}\r@gistered
                                                \else$\m@th\r@gistered$\fi}%
 \def\r@gistered{^{\ooalign
  {\hfil\raise.07ex\hbox{$\scriptstyle\rm\text{R}$}\hfil\crcr
  \mathhexbox20D}}}}{}%
\@ifundefined{Eth}{}{}%
\@ifundefined{eth}{}{}%
\@ifundefined{Thorn}{}{}%
\@ifundefined{thorn}{}{}%
\@ifundefined{degree}{}{}%
\def\BibTeX{{\rm B\kern-.05em{\sc i\kern-.025em b}\kern-.08em
    T\kern-.1667em\lower.7ex\hbox{E}\kern-.125emX}}%
\newdimen\theight
\def\Column{%
 \vadjust{\setbox\z@=\hbox{\scriptsize\quad\quad tcol}%
  \theight=\ht\z@\advance\theight by \dp\z@\advance\theight by \lineskip
  \kern -\theight \vbox to \theight{%
   \rightline{\rlap{\box\z@}}%
   \vss
   }%
  }%
 }%
\def\qed{%
 \ifhmode\unskip\nobreak\fi\ifmmode\ifinner\else\hskip5\p@\fi\fi
 \hbox{\hskip5\p@\vrule width4\p@ height6\p@ depth1.5\p@\hskip\p@}%
 }%
\def\miss{\hbox{\vrule height2\p@ width 2\p@ depth\z@}}%
\def\tcol#1{{\baselineskip=6\p@ \vcenter{#1}} \Column}  %
\makeatother
\newtheorem{theorem}{Theorem}
\newtheorem{lemma}[theorem]{Lemma}

\newtheorem{definition}[theorem]{Definition}

\begin{document}

\author{Paul C. Eklof}
\address{Department of Mathematics, University of California, Irvine, 
Irvine, CA 92717}
\email{pceklof@@uci.edu}
\thanks{Thanks to Rutgers University for 
funding  the authors' visits to Rutgers.}
\author{Saharon Shelah}
\thanks{The second author was partially supported by the {\it Basic Research 
Foundation} administered by The Israel Academy of Sciences and 
Humanities. Pub. No. 505. Final version.}
\address{Institute of Mathematics, Hebrew University, Jerusalem 91904, 
Israel}
\email{shelah@@math.huji.ac.il}

\title[A Combinatorial Principle]{A Combinatorial Principle  Equivalent to the Existence of Non-free Whitehead
Groups }
\date{}


\maketitle
\begin{abstract}
As a consequence of identifying the principle described in the title, we
prove that for any uncountable cardinal $\lambda $, if there is a $\lambda $%
-free Whitehead group of cardinality $\lambda $ which is not free, then
there are many ``nice'' Whitehead groups of cardinality $\lambda $ which are
not free.
\end{abstract}


\section{Introduction}

Throughout, ``group'' will mean abelian group; in particular, ``free group''
will mean free abelian group.

Two problems which have been shown to be undecidable in ZFC (ordinary set
theory) for some uncountable $\lambda $ are the following:

\begin{itemize}
\item  Is there a group of cardinality $\lambda $ which is $\lambda $-free
(that is, every subgroup of cardinality $<\lambda $ is free), but is not
free?

\item  Is there a Whitehead group $G$ (that is, $\limfunc{Ext}(G,{\Bbb Z})=0$%
) of cardinality $\lambda $ which is not free?
\end{itemize}

\noindent (See \cite{MS} for the first, \cite{Sh74} for the second; also \cite{EM} is
a general reference for unexplained terminology and further information.)

The second author has shown that the first problem is equivalent to a
problem in pure combinatorial set theory (involving the 
important notion of a $\lambda$-system; see Theorem \ref{trans}.) This
not only makes it easier to prove independence results (as in \cite{MS}),
but also allows one to prove (in ZFC!) group-theoretic results such as:
\smallskip

\begin{quote}
($\bigtriangledown $) if there is a $\lambda $-free group of cardinality $%
\lambda $ which is not free, then there is a strongly $\lambda $-free group
of cardinality $\lambda $ which is not free.
\end{quote}

\noindent (See \cite{Sh85} or \cite[Chap. VII]{EM}.) A group $G$ is said to
be {\it strongly $\lambda $-free} if every subset of $G$ of cardinality $%
<\lambda $ is contained in a free  subgroup $H$ of cardinality $<\lambda $
such that $G/H$ is $\lambda $-free. One reason for interest in this class of
groups is that they are precisely the groups which are equivalent to a free
group with respect to the infinitary language $L_{\infty \lambda }$ (see 
\cite{E}).

There is no known way to prove ($\bigtriangledown $) except to go through
the combinatorial equivalent.

As for the second problem, the second author has shown that for $\lambda
=\aleph _1$, there is a combinatorial characterization of the problem:

\begin{quote}
there is a non-free Whitehead group of cardinality $\aleph _1$ if and only
if there is a ladder system on a stationary subset of $\aleph _1$ which
satisfies 2-uniformization.
\end{quote}

\noindent (See the Appendix to this paper; a knowledge of the undefined
terminology in this characterization is not needed for the body of this
paper.) Again, there are group-theoretic consequences which are provable in
ZFC:
\smallskip

\begin{quote}
($\bigtriangledown \bigtriangledown $) if there is a non-free Whitehead
group of cardinality $\aleph _1$, then there is a strongly $\aleph _1$-free,
non-free Whitehead group of cardinality $\aleph _1$.
\end{quote}
\smallskip

\noindent (See \cite[\S XII.3]{EM}. It is consistent with ZFC that there are
Whitehead groups of cardinality $\aleph _1$ which are not strongly $\aleph
_1 $-free.)

Our aim in this paper is to generalize ($\bigtriangledown \bigtriangledown $%
) to cardinals $\lambda >\aleph _1$ by combining the two methods used to
prove ($\bigtriangledown $) and ($\bigtriangledown \bigtriangledown $).
Since the existence of a non-free Whitehead group $G$ of cardinality $\aleph
_1$ implies that for every uncountable cardinal $\lambda $ there exist
non-free Whitehead groups of cardinality $\lambda $ (e.g. the direct sum of $%
\lambda $ copies of $G$ --- which is not $\lambda $-free), the appropriate
hypothesis to consider is:

\begin{itemize}
\item  there is a $\lambda $-free Whitehead group of cardinality $\lambda $
which is not free.
\end{itemize}

By the Singular Compactness Theorem (see \cite{Sh75}), $\lambda $ must be
regular. It can be proved consistent that there are uncountable $\lambda
>\aleph _1$ such that the hypothesis holds. (See \cite{ES} or \cite{T}.) In
particular, for many $\lambda $ (for example, $\lambda =\aleph _{n+1}$, $%
n\in \omega $) it can be proved consistent with ZFC that $\lambda $ is the
smallest cardinality of a non-free Whitehead group; hence there is a $%
\lambda $-free Whitehead group of cardinality $\lambda $ which is not free.
(See \cite{ES}.)

Our main theorem is then:

\begin{theorem}
\label{main}If there is a $\lambda $-free Whitehead group of cardinality $%
\lambda $ which is not free, then there are $2^\lambda $ different strongly $%
\lambda $-free Whitehead groups of cardinality $\lambda $.
\end{theorem}

Our proof will proceed in three steps. First, assuming the hypothesis ---
call it (A) --- of the Theorem, we will prove

\begin{quote}
(B) there is a combinatorial object, consisting of a $\lambda $-system with
a type of uniformization property.
\end{quote}

Second, we will show that the combinatorial property (B) can be improved to
a stronger combinatorial property (B+), which includes the ``reshuffling
property''. Finally, we prove that (B+) implies 
$$
\mbox{(A+): the existence of many strongly }\lambda \mbox{-free Whitehead
groups.} 
$$

Note that we have found, in (B) or (B+), a combinatorial property which is
equivalent to the existence of a non-free $\lambda $-free Whitehead group of
cardinality $\lambda $, answering an open problem in \cite[p. 453]{EM}.
Certainly this combinatorial characterization is more complicated than the
one for groups of cardinality $\aleph _1$ cited above. This is not
unexpected; indeed the criterion for the existence of $\lambda $-free groups
in Theorem \ref{trans} implies that the solution to the open problem is
inevitably going to involve the notion of a $\lambda $-system. A good reason
for asserting the interest of the solution is that it makes possible the
proof of Theorem \ref{main}.

In an Appendix we provide a  simpler proof than the previously published one
of the fact that the existence of a non-free Whitehead group of cardinality $%
\aleph _1$ implies the existence of a ladder system on a stationary subset
of $\aleph _1$ which satisfies $2$-uniformization.

We thank Michael O'Leary for his careful reading of and comments on this paper.

\section{Preliminaries}

The following notion, of a $\lambda $-set, may be regarded as a
generalization of the notion of a stationary set.

\begin{definition}

   (1) The set of  all functions $$\eta \colon n =
 \{0,\ldots, n-1\} \rightarrow  \lambda$$  $(n \in  \omega )$ is denoted 
$^{<\omega }\!\lambda $; the domain of $\eta $ is denoted $\ell (\eta )$ 
and called the {\em length} of $\eta $; we identify $\eta $ with the sequence
$$
\langle \eta (0), \eta (1),\ldots, \eta (n-1)\rangle .
$$
 
\noindent
Define a partial ordering  on $^{<\omega }\!\lambda $ by:  $\eta _1 \leq  
\eta _2$ if and only if $\eta _1$ is a restriction of $\eta _2$. This makes 
$^{<\omega }\!\lambda $ into a {\em tree}. For any  $\eta  = 
\langle \alpha _0, \ldots, \alpha _{n-1}\rangle  \in {}^{<\omega }\!\lambda $, 
$\eta \smallfrown\langle \beta \rangle $ denotes the sequence 
$\langle \alpha _0,\ldots, \alpha _{n-1}, \beta \rangle $. If $S$ is a 
subtree of $^{<\omega }\!\lambda $, an element $\eta $ of $S$ is called a 
{\em final node} of $S$ if no $\eta \smallfrown\langle \beta \rangle $ belongs to 
$S$. Denote the set of final nodes of $S$ by $S_f$.

\medskip
 
(2) A $\lambda$-{\em set} is a subtree $S$ of $^{<\omega }\!\lambda $ together with a 
cardinal $\lambda _\eta $ for every  $\eta  \in  S$  such that  
$\lambda _\emptyset  = \lambda $,  and: 
 
(a) for all  $\eta  \in  S$, $\eta $ is a final node of $S$ if and only if  
$\lambda _\eta  = \aleph _0;$ 
 
(b) if $\eta \in S \setminus S_f{ }$, then $\eta \smallfrown\langle \beta
\rangle \in S$ implies $\beta \in \lambda _\eta $, $\lambda _{\eta
\smallfrown\langle \beta \rangle } < \lambda _\eta $ and $E_\eta  =  \{\beta
< \lambda _\eta \colon \eta \smallfrown\langle \beta \rangle \in S\}$ is
stationary in $\lambda _\eta $. 
 
\medskip
 
(3) A $\lambda ${\em -system} is a $\lambda$-set together with a set $B_\eta $
for each $\eta  \in  S$ such that $B_\emptyset  = \emptyset $, and for all 
$\eta  \in  S \setminus  S_f$: 
 
(a) for all  $\beta  \in  E_\eta $,  
$\lambda _{\eta \smallfrown\langle \beta \rangle } \leq  
|B_{\eta \smallfrown\langle \beta \rangle }| < \lambda _\eta ;$ 
 
(b) $\{B_{\eta \smallfrown\langle \beta \rangle }\colon \beta \in E_\eta \}$
is a continuous chain of sets, i.e.\  if $\beta \leq \beta '$ are in $E_\eta
$, then $B_{\eta \smallfrown\langle \beta \rangle } \subseteq B_{\eta
\smallfrown\langle \beta '\rangle }$, and if $\sigma $ is a limit point of
$E_\eta $, then $B_{\eta \smallfrown\langle \sigma \rangle } = \cup \{B_{\eta
\smallfrown\langle \beta \rangle }\colon \beta < \sigma $, $\beta \in E_\eta
\}$;

 \medskip 
 
(4) For any $\lambda $-system $\Lambda = (S$, $\lambda _\eta $, $B_\eta
\colon \eta \in S)$, and any $\eta \in S$, let $\bar B_\eta = \cup
\{B_{\eta \rest m}\colon m \leq \ell (\eta )\}$.  Say that a family
${\cal S} = \{s_\zeta \colon \zeta \in S_f\}$ of countable sets is {\em based on $\Lambda $} if ${\cal S}$
is indexed by $S_f{ }$ and for every $\zeta \in S_f{ }$, $s_
\zeta
\subseteq \bar B_\zeta $.

\medskip

(5) A family ${\cal S}=\{s_i:i\in I\}$ is said to be {\em free} if it has a
transversal, that is, a one-one function $T:I\rightarrow \cup {\cal S}$ such
that for all $i\in I$, $T(i)\in s_i$. We say ${\cal S\ }$is $\lambda ${\em %
-free} if every subset of ${\cal S}$ of cardinality $<\lambda $ has a
transversal.
\end{definition}

It can be proved that a family ${\cal S}$ which is based on a $\lambda $%
-system is not free. (See \cite[Lemma VII.3.6]{EM}.)

\smallskip\ \ 

The following theorem now gives combinatorial equivalents to the existence
of a $\lambda $-free group of cardinality $\lambda $ which is not free. (See 
\cite{Sh85} or \cite[\S VII.3]{EM}.)

\begin{theorem}
\label{trans} For any uncountable cardinal $\lambda $, the following are
equivalent:

\begin{enumerate}
\item  there is a $\lambda $-free group of cardinality $\lambda $ which is
not free;

\item  there is a family ${\cal S}$ of countable sets such that ${\cal S}$
has cardinality $\lambda $ and is $\lambda $-free but not free;

\item  there is a family ${\cal S}$ of countable sets such that ${\cal S}$
has cardinality $\lambda $, is $\lambda $-free, and is based on a $\lambda $%
-system.
\end{enumerate}
\end{theorem}

\begin{definition}

 A subtree $S$ of $^{<\omega }\!\lambda $ is said to have {\em height
$n$} if all the final nodes of $S$ have length $n$.  A $\lambda$-set or
$\lambda $-system is said to have height $n$ if its associated subtree
$S$ has height $n$.

\end{definition}

A $\lambda $-set of height 1 is essentially just a stationary subset of $%
\lambda $. Not every $\lambda $-set has a height, but the following lemma
implies that every $\lambda $-set contains one which has a height. It is a
generalization, and a consequence, of the fact that if a stationary set $E$
is the union of subsets $E_n$ ($n\in \omega $), then for some $n$, $E_n$ is
stationary (cf. \cite[Exer. 2, p. 238]{EM}).

\begin{lemma}
\label{encap2}

 If $(S, \lambda _\eta : \eta \in S)$ is a 
$\lambda $-set, and  $S_f = \bigcup _{n\in \omega } S^n_f{ }$, let  $S^n = 
\{\eta  \in  S\colon  \eta  \leq  \tau $ for some  $\tau  \in  S^n_f{ }\}$.
Then 
 for some $n$, $(S^n, \lambda _\eta : \eta \in S^n)$
contains a $\lambda $-set.
\end{lemma}

If $\Lambda =(S,\lambda _\eta ,B_\eta :\eta \in S)$ is a $\lambda $-system
of height $n$, and ${\cal S}=\{s_\zeta :\zeta \in S_f\}$ is a family of
countable sets based on $\Lambda $, let $s_\zeta ^k=s_\zeta \cap B_{\zeta
\rest k} $ for $1\leq k\leq n$.

\smallskip\ 

The following is useful in carrying out an induction on $\lambda $-systems.

\begin{definition}
\label{genlsys}

 Given a $\lambda $-system $\Lambda = (S, \lambda _\eta , B_\eta
\colon \eta \in S)$ and a node $\eta $ of $S$,
 let $S^\eta = \{\nu \in S\colon \eta \leq \nu \}$.
We will denote by $\Lambda^\eta$ the   $\lambda _\eta $-system 
 which is naturally isomorphic to
$(S^\eta , \lambda _\nu , B'_\nu \colon \nu \in S^\eta )$
where $B'_\eta = \emptyset $ and $B'_\nu = B_\nu $ if $\nu \neq
\eta$.

{\em (That is, we replace the initial node, $\eta$, of $S^\eta$ by
$\emptyset$, and translate the other nodes accordingly.)}

If ${\cal S} = \{s_\zeta \colon \zeta \in S_f\}$ is a family of
 countable sets  {\em based on} $\Lambda $ and $\zeta \in
S^\eta_f$,
 let $s^\eta _\zeta = \cup
\{s^k_\zeta \colon k > \ell (\eta )\}$.  Let ${\cal S}^\eta = \{s^\eta
_\zeta \colon \zeta \in S^\eta_f \}$; it is a family of countable sets based
on $\Lambda ^\eta .$

\end{definition}

In order to construct a (strongly) $\lambda $-free group from a family of
countable sets based on a $\lambda $-system, we need that the family have an
additional property:

\begin{definition}
\label{reshuff}A family ${\cal S}$ of countable sets based on a $\lambda $%
-system $\Lambda $ is said to have the {\em {reshuffling property}} if for
every $\alpha <\lambda $ and every subset $I$ of $S_f$ such that $%
|I|<\lambda $, there is a well-ordering $<_I$ of $I$ such that for every $%
\tau $, $\zeta $$\in I$, $s_\zeta \setminus \bigcup_{\nu <_I\zeta }s_\nu $
is infinite, and $\tau (0)\leq \alpha <\zeta (0)$ implies that $\tau
<_I\zeta .$
\end{definition}

It can be shown (in fact it is part of the proof of the theorem) that the
three equivalent conditions in Theorem \ref{trans} are equivalent to:

\begin{quote}
(iv) {\it there is a family ${\cal S}$ of countable sets such that ${\cal S}$
has cardinality $\lambda $, is $\lambda $-free, is based on a $\lambda $%
-system, and${\cal \ }$ has the reshuffling property.}
\end{quote}

\ 

Finally, for future reference, we observe the following simple fact:

\begin{lemma}
\label{nonfree}Suppose that for some integers $r\geq 0$, and $%
d_m^\ell $, and some primes $q_m$ ($\ell <r$, $m\in \omega $), $H$ is the
abelian group on the generators $\{z_j:j\in \omega \}$ modulo the relations 
$$
q_mz_{m+r+1}=z_{m+r}+\sum_{\ell <r}d_m^\ell z_\ell  
$$
for all $m\in \omega $. Then $H$ is not free.

Conversely, if $C$ is a torsion-free abelian group of rank $r+1$ which is not free but
is such that every subgroup of rank $\leq r$ is free, then $C$ contains a
subgroup $H$ which is given by generators and relations as above.
\end{lemma}

{\sc proof. }Let{\sc \ }$H$ be as described in the first part. If{\sc \ }$H$
is free, then $H$ is finitely generated, since it clearly has rank $\leq r+1$%
. {\sc \ }Let $L$ be the subgroup of $H$ generated by (the images of) $%
z_0,...,z_{r-1}$. By comparing coefficients of linear combinations in the
free group on $\{z_j:j\in \omega \}$, one can easily verify that $L$ is a
pure subgroup of $H$, and that $H/L$ is a rank one group which is not free
(because $z_r+L$ is non-zero and divisible by $q_0q_1\cdot \cdots \cdot q_m$
for all $m\in \omega $) and hence not finitely-generated. But this is
impossible if $H$ is free.

Conversely, let $C$ be as stated, and let $L$ be a pure subgroup of rank $r$%
. Then $L$ is free (say with basis $z_0,...,z_{r-1}$) and $C/L$ is a
non-free torsion-free group of rank $1$. Thus $C/L$ contains a subgroup with
a non-zero element $z_r+L$ such that either: $z_r+L$ is divisible by all
powers of $p$ for some prime $p$ (in which case we let $q_m=p$ for all $m$);
or $z_r+L$ is divisible by infinitely many primes (in which case we let $%
\{q_m:m\in \omega \}$ be an infinite set of primes dividing $z_r$). It is
then easy to see that $H$ exists as desired. $\Box $



\section{(A) implies (B)}

\begin{theorem}
\label{AtoB}

For any regular uncountable cardinal $\lambda $, if

(A) there is a Whitehead group of cardinality $\lambda $ which is $\lambda $%
-free but not free,

then

(B) there exist integers $n>0$ and $r\geq 0$, and:

\begin{enumerate}
\item  a $\lambda $-system $\Lambda =(S,\lambda _\eta ,B_\eta :\eta \in S)$
of height $n$;

\item  one-one functions $\varphi _\zeta ^k$ $(\zeta \in S_f$, $1\leq k\leq
n)$ with $\limfunc{dom}(\varphi _\zeta ^k)=\omega $;

\item  primes $q_{\zeta ,m}$ ($\zeta \in S_f$, $m\in \omega $); and

\item  integers $d_{\zeta ,m}^\ell $ ($\zeta \in S_f$, $m\in \omega $, $\ell
<r$)
\end{enumerate}

\noindent 

such that

\begin{quote}
(a) if we define $s_\zeta =\bigcup_{k=1}^n\limfunc{rge}(\varphi _\zeta ^k)$,
then ${\cal \ S}=\{s_\zeta :\zeta \in S_f\}$ is a $\lambda $-free family of
countable sets based on $\Lambda $; in particular, $\limfunc{rge}(\varphi
_\zeta ^k)\subseteq B_{\zeta \rest k}$;
\end{quote}

\noindent 

and

\begin{quote}
(b) for any functions $c_\zeta :\omega \to {\Bbb Z}$ ($\zeta \in S_f$), there
is a function $f:\bigcup {\cal S}\to {\Bbb Z\ }$such that for all $\zeta \in
S_f$ there are integers $a_{\zeta ,j}$ ($j\in \omega $) such that for all $%
m\in \omega $, 
$$
c_\zeta (m)=q_{\zeta ,m}a_{\zeta ,m+r+1}-a_{\zeta ,m+r}-\sum_{\ell
<r}d_{\zeta ,m}^\ell a_{\zeta ,\ell }-\sum_{k=1}^nf(\varphi _\zeta ^k(m)). 
$$
\end{quote}
\end{theorem}

{\sc proof. }We shall refer to the data in (B), which satisfies (a) and (b),
as a {\it \ $\lambda $-system with data for the Whitehead problem} or more
briefly a {\it Whitehead $\lambda $-system}. Given a Whitehead group $G$ of
cardinality $\lambda $ which is $\lambda $-free but not free,{\sc \ }we
begin by defining a $\lambda $-system and a family of countable sets based
on the $\lambda $-system following the procedure given in \cite[VII.3.4]{EM}%
; we review that procedure here.

Choose a {\em $\lambda$-filtration} of $G$, that is,  write $G$ as
the union of a continuous chain $$
G = \bigcup _{\alpha <\lambda } B_\alpha 
$$
 
\noindent
of pure subgroups of cardinality $< \lambda$ such that if  
$G/B_\alpha $ is not $\lambda $-free, then $B_{\alpha +1}/B_\alpha $ is not 
free. Since $G$ is not free,
 $$ E_\emptyset  = \{\alpha  < \lambda
\colon B_\alpha \hbox{ is  not }\ \lambda \hbox{-pure in }\ G\}
$$
 
\noindent
is stationary in $\lambda$. For each $\alpha  \in  E_\emptyset $,
let $\lambda _\alpha$ $(< \lambda )$ be minimal such that $B_{\alpha
+1}/B_\alpha $ has a subgroup of cardinality $\lambda _\alpha $ which
is not free; $\lambda _\alpha $ is regular by the Singular
Compactness Theorem (see \cite{Sh75} or \cite[IV.3.5]{EM}). 
 If $\lambda _\alpha $ is countable, then let $\langle \alpha
\rangle$ 
be a final node of the tree; otherwise choose $G_\alpha \subseteq
B_{\alpha +1}$ of cardinality $\lambda _\alpha $ such that
 
$$H_\alpha =  (G_\alpha  + B_\alpha )/B_\alpha $$
is not free. Then $H_\alpha $ is $\lambda _\alpha $-free, and we can choose a 
$\lambda _\alpha $-filtration of $G_\alpha ,$
$$
G_\alpha  = \bigcup _{\beta <\lambda _\alpha } B_{\alpha ,\beta } $$
 
\noindent
such that for all $\beta$, $(B_{\alpha ,\beta } + B_\alpha )/B_\alpha $
is pure in $H_\alpha$ and  if $(B_{\alpha ,\beta } + B_\alpha )/B_\alpha $ is not 
$\lambda _\alpha $-pure in $H_\alpha $, then 
$$
(B_{\alpha ,\beta +1} + B_\alpha /B_\alpha )/(B_{\alpha ,\beta } + 
B_\alpha /B_\alpha ) \cong  (B_{\alpha ,\beta +1} + 
B_\alpha )/(B_{\alpha ,\beta } + B_\alpha )
$$
 
\noindent
is not free. Since $H_\alpha $ is not free,
$$
E_\alpha = \{\beta  < \lambda _\alpha \colon B_{\alpha ,\beta } + 
B_\alpha /B_\alpha \hbox{ is not }\ \lambda _\alpha \hbox{-pure 
in }\ H_\alpha \}
$$
 
\noindent
is stationary in $\lambda _\alpha $. For each  $\beta  \in  E_\alpha $, choose 
$\lambda _{\alpha ,\beta }\ (< \lambda _\alpha )$ minimal such that there is a 
subgroup $G_{\alpha ,\beta }$ of $B_{\alpha ,\beta +1}$ of cardinality 
$\lambda _{\alpha ,\beta }$ so that
$$
H_{\alpha ,\beta } = (G_{\alpha ,\beta } + B_{\alpha ,\beta } + 
B_\alpha )/(B_{\alpha ,\beta } + B_\alpha )
$$
 
\noindent
is not free. If $\lambda _{\alpha ,\beta }$ is countable, let
$\langle \alpha, \beta \rangle $ be a final node;  otherwise choose
a $\lambda _{\alpha ,\beta }$-filtration of $G_{\alpha ,
\beta }$.
Continue in this way along each branch until a final node is reached.

 As we have just done,
we will use, when convenient, the notation $G_{\alpha ,\beta }$ instead of $%
G_{\langle \alpha ,\beta \rangle }$, etc.; thus for example we
will write $G_{\eta ,\delta }$ instead of $G_{\eta \frown \langle
\delta \rangle }$.

In this way we obtain a $\lambda $-system $\Lambda =(S,\lambda _\eta ,B_\eta
:\eta \in S)$ where for each $\zeta $$\in S_f$, there is a countable
subgroup $G_\zeta $ of $G$ such that 
$$
G_\zeta +\langle \bar B_\zeta \rangle /\langle \bar B_\zeta \rangle 
$$
is not free. We can assume that for each $\eta \in S\setminus S_f$ and each $%
\delta \in E_\eta $, 
$$
B_{\eta ,\delta +1}+\langle \bar B_\eta \rangle =G_{\eta
,\delta }+B_{\eta ,\delta }+\langle \bar B_\eta \rangle 
\mbox{.}
$$

\noindent
We can also assume  that for all
$\zeta \in S_f$, $G_\zeta$ has been chosen so that
 $G_\zeta  + \langle \bar B_\zeta \rangle /\langle \bar B_\zeta \rangle$ has
  finite rank $r_\zeta + 1$ for some $r_\zeta$ such that every subgroup
 of rank $\leq r_\zeta $ is free. By restricting to a sub-$\lambda$-set,
 we can assume that there is an $r$ such that $r_\zeta +1 = r + 1$
for all 
  $\zeta \in S_f$ and that there is an $n$ such that $\Lambda$ has
height $n$ (cf. Lemma \ref{encap2}). Moreover, we can assume (easing the 
purity condition, if necessary) that
$G_\zeta + \langle \bar B_\zeta \rangle /\langle \bar B_\zeta \rangle$
is as described in  Lemma \ref{nonfree}, that is,
 it is generated
modulo $\langle \bar B_\zeta \rangle$ by the cosets
 of elements $z_{\zeta, j}$ which satisfy precisely the relations
which are consequences of relations  
 $$q_{\zeta ,m}z_{\zeta ,m+r+1}=z_{\zeta ,m+r}+\sum_{\ell <r}d_{\zeta
,m}^\ell z_{\zeta ,\ell }$$ (modulo $\langle \bar B_\zeta
\rangle$) for some primes $q_{\zeta,m}$ and integers $d^\ell_{\zeta
,m}$. Fix $g_{\zeta,m} \in \langle \bar B_\zeta \rangle$ such that
in $G$ $$q_{\zeta ,m}z_{\zeta ,m+r+1}=z_{\zeta ,m+r}+\sum_{\ell
<r}d_{\zeta ,m}^\ell z_{\zeta ,\ell }+ g_{\zeta,m}.$$

There is a countable subset $t_\zeta $ of $\bar B_\zeta $ such
that $G_\zeta  \cap  \langle \bar B_\zeta \rangle $ is contained in the subgroup
generated by $t_\zeta $. Let  $s_\zeta  = t_\zeta  \times  \omega $. Then it is proved in
\cite[VII.3.7]{EM} that $\{s_\zeta \colon \zeta  \in  S_f{ }\}$  is
$\lambda $-free and based on the $\lambda$-system 
$ (S,\lambda _\eta
,B_\eta ^{\prime }:\eta \in S)$ where $B_\eta ^{\prime }=B_\eta
\times \omega $. 

Let $s_\zeta ^k=s_\zeta \cap B_{\zeta \rest k}^{\prime }$ and let $\nu _\zeta
^k:\omega \rightarrow $ $s_\zeta ^k$ enumerate $s_\zeta ^k$ without
repetition. We can write each $g_{\zeta ,m}$ as a sum $\sum_{k=1}^ng_{\zeta
,m}^k$ where $g_{\zeta ,m}^k\in B_{\zeta \rest k}^{}$. Now define
$$
\varphi _\zeta ^k(m)=\langle \nu _\zeta ^k(m),g_{\zeta
,m}^k\rangle \in B_{\zeta \rest k}^{\prime }\times B_{\zeta \rest k}\mbox{.}
$$
Then $s_\zeta ^{\prime \prime }=$ $\bigcup_{k=1}^n\limfunc{rge}(\varphi
_\zeta ^k)$ is based on the $\lambda $-system $(S,\lambda _\eta ,B_\eta
^{\prime \prime }:\eta \in S)$ where $B_\eta ^{\prime \prime }=B_\eta
^{\prime }\times B_\eta $. Moreover $\{s_\zeta ^{\prime \prime }:\zeta \in
S_f\}$ is a $\lambda $-free family because of the choice of the first
coordinate of $\varphi _\zeta ^k(m)$. Thus we have defined the data in (B)
such that (a) holds. It remains to verify (b). So let $c_\zeta :\omega \to 
{\Bbb Z}$ ($\zeta \in S_f$) be given. We are going to define a short exact
sequence%
$$
0\longrightarrow {\Bbb Z}\longrightarrow M\stackrel{\pi }{\longrightarrow }%
G\longrightarrow 0 
$$
and then use a splitting of $\pi $ to define the function $f:\bigcup {\cal S}%
\to {\Bbb Z}$.

We will use the lexicographical ordering, $<_\ell $, on $S$ defined as
follows: $\eta _1<_\ell \eta _2$ if and only if either $\eta _1$ is a
restriction of $\eta _2$ or $\eta _1(i)<\eta _2(i)$ for the least $i$ such
that $\eta _1(i)\neq \eta _2(i)$. Note that if $\eta _1<_\ell \eta _2$, then 
$\langle \bar B_{\eta _1}\rangle $ $\subseteq \langle \bar B%
_{\eta _2}\rangle $. The lexicographical ordering is a well-ordering of $S$,
so there is an order-preserving bijection $\theta :\tau \rightarrow
\langle S,<_\ell \rangle $ for some ordinal $\tau $. If for each $%
\sigma <\tau $ we let $A_\sigma =\langle \bar B_{\theta (\sigma
)}\rangle $, then $G=\cup _{\sigma <\tau }A_\sigma $ represents $G$ as the
union of a chain of subgroups. However, we must exercise caution since, as
we will see, this chain is not necessarily continuous.

The kernel of $\pi $ will be generated by an element $e\in M$. We will
define $\pi $ to be the union of a chain of homomorphisms $\pi _\sigma
:M_\sigma \rightarrow A_\sigma \rightarrow 0$ with kernel ${\Bbb Z}e$. The $%
\pi _\sigma $ will be defined by induction on $\sigma $. At the same time,
we will also define, as we go along, a chain of set functions $\psi _\sigma
:A_\sigma \rightarrow M_\sigma $ such that $\pi _\sigma \circ \psi _\sigma =%
\limfunc{id}_{A_\sigma }$. Let $\pi _0$ be the zero homomorphism $:{\Bbb Z}%
e\rightarrow A_0=\{0\}$.

Suppose that $\pi _\rho $ and $\psi _\rho $ have been defined for all $\rho
<\sigma $ for some $\sigma <\tau $; say $\theta (\sigma )=\eta $ where $\eta 
$$=$$\langle \nu ,\delta \rangle $ for some $\nu \in S$, $\delta \in E_\nu $%
. Suppose first that $\sigma $ is a limit ordinal. Let $\pi _\sigma ^{\prime
}:M_\sigma ^{\prime }\rightarrow \cup _{\rho <\sigma }A_\rho $ be the direct
limit of the $\pi _\rho $ ($\rho <\sigma )$ and let $\psi _\sigma ^{\prime }$
be the direct limit of the $\psi _\rho $. In particular, $M_\sigma ^{\prime
}=\underrightarrow{\lim }\{M_\rho :\rho <\sigma \}$. If $\cup _{\rho <\sigma
}A_\rho =A_\sigma $, then we can let $\pi _\sigma =\pi _\sigma ^{\prime }$
and $\psi _\sigma =\psi _\sigma ^{\prime }$; this will happen, for example,
if $\delta $ is a limit point of $E_\nu $.

But it may be that $\delta $ has an immediate predecessor $\delta _1\in
E_\nu $. (Since $\sigma $ is a limit ordinal, it follows that $\eta $ is not
a final node of $S$.) Then 
$$
\cup _{\rho <\sigma }A_\rho =\cup _{\gamma <\lambda _{\nu ,\delta _1}}B_{\nu
,\delta _1,\gamma }+B_{\nu ,\delta _1}+\langle \bar B_\nu
\rangle =B_{\nu ,\delta _1+1}+\langle \bar B_\nu
\rangle \mbox{.}
$$
Notice that $\cup _{\rho <\sigma }A_\rho $ will be a proper subgroup of $%
A_\sigma $ if $\delta _1+1<\delta $ (i.e. if $\delta _1+1\notin E_\nu $). We
can extend $\pi _\sigma ^{\prime }$ to $\pi _\sigma :M_\sigma \rightarrow
A_\sigma $ because the inclusion of $\cup _{\rho <\sigma }A_\rho $ into $%
A_\sigma $ induces a surjection of $\limfunc{Ext}(A_\sigma ,{\Bbb Z)}$ onto $%
\limfunc{Ext}(\cup _{\rho <\sigma }A_\rho ,{\Bbb Z)}$. Finally, extend $\psi
_\sigma ^{\prime }$ to $\psi _\sigma $ in any way such that $\pi _\sigma
\circ \psi _\sigma $ is the identity on $A_\sigma $.

Now let us consider the case when $\sigma =\rho +1$ is a successor ordinal.
Recall that $\theta (\sigma )=\langle \nu ,\delta \rangle $. There are two
subcases. In the first, $\delta $ is the least element of $E_\nu $, so $%
\theta (\rho )=$ $\nu $ and $A_\rho =\langle \bar B_\nu
\rangle $, $A_\sigma =B_{\nu ,\delta }+\langle \bar B_\nu
\rangle $. In this subcase, we extend $\pi _\rho $ to $\pi _\sigma $
using the surjectivity of $\limfunc{Ext}(A_\rho ,{\Bbb Z)}\rightarrow $ $%
\limfunc{Ext}(A_\sigma ,{\Bbb Z)}$.

In the second and last subcase, $\delta $ has an immediate predecessor $%
\delta _1$ in $E_\nu $; then $\theta (\rho )=$$\langle v,\delta
_1\rangle $, a final node of $S$. Let $\zeta $ denote $\langle
v,\delta _1\rangle $; then $B_{\nu ,\delta _1+1}+\langle \bar B_\zeta
\rangle /\langle \bar B_\zeta \rangle $ is as described in Lemma \ref
{nonfree}, that is, it is generated modulo $\langle \bar B_\zeta \rangle $
by the cosets of elements $z_{\zeta ,j}$ which satisfy the relations 
\begin{align}
\label{zzeta}
q_{\zeta ,m}z_{\zeta ,m+r+1}=z_{\zeta 
,m+r}+\sum_{\ell<r}
d_{\zeta ,m}^\ell z_{\zeta ,\ell }+\sum_{k=1}^ng_{\zeta ,m}^k 
\end{align}
in $G$ for some primes $q_{\zeta ,m}$, integers $d_{\zeta ,m}^\ell $ and
elements $g_{\zeta ,m}^k\in B_{\zeta \rest k}$. It is at this point that we use
the function $c_\zeta $. Define $M_\sigma ^{\prime }$ to be generated over $%
M_\rho $ by elements $z_{\zeta ,j}^{\prime }$ modulo the relations 
\begin{align}
\label{zprime}q_{\zeta ,m}z_{\zeta ,m+r+1}^{\prime }=z_{\zeta ,m+r}^{\prime
}+\sum_{\ell <r}d_{\zeta ,m}^\ell z_{\zeta ,\ell }^{\prime
}+\sum_{k=1}^n\psi _\rho (g_{\zeta ,m}^k)+c_\zeta (m)e 
\end{align}
and define 
$$
\pi _\sigma ^{\prime }:M_\sigma ^{\prime }\rightarrow B_{\nu ,\delta
_1+1}+\langle \bar B_\zeta \rangle 
$$
to be the homomorphism extending $\pi _\rho $ which takes $z_{\zeta
,j}^{\prime }$ to $z_{\zeta ,j}$. One can verify that $\pi _\sigma ^{\prime
} $ is well-defined and has kernel ${\Bbb Z}e$. Extend $\psi _\rho $ to $\psi
_\sigma ^{\prime }$ in any way such that $\pi _\sigma ^{\prime }\circ \psi
_\sigma ^{\prime }$ is the identity. We extend $\pi _\sigma ^{\prime }$ to $%
\pi _\sigma :M_\sigma \rightarrow A_\sigma =\langle \bar B%
_{\langle \nu ,\delta \rangle }\rangle $ by using the
surjectivity of $\limfunc{Ext}(A_\sigma ,{\Bbb Z})\rightarrow $ $\limfunc{Ext}%
(B_{\nu ,\delta _1+1}+\langle \bar B_\zeta \rangle ,{\Bbb Z}%
) $; finally we extend $\psi _\sigma ^{\prime }$.

This completes the definition of $\pi :M\rightarrow G$ and of the set map $%
\psi :G\rightarrow M$ ($=$ the direct limit of the $\psi _\sigma $). Since $%
G $ is a Whitehead group, there is a homomorphism $\rho $$:G\rightarrow M$
such that $\pi \circ \rho $ is the identity on $G$. In order to define $f$,
consider an element $x$ of $\cup {\cal S}$; $x$ is an ordered pair equal to $%
\varphi _\zeta ^k(m)$ (possibly for many different $(\zeta ,k,m)$). If $g$
is the second coordinate of $x$, let $f(x)$ be the unique integer such that 
$$
\psi (g)-\rho (g)=f(x)e\mbox{.}
$$
Also for any $\zeta $$\in S_{f\mbox{ }}$ and $j\in \omega $ define $a_{\zeta
,j}$ such that 
$$
z_{\zeta ,j}^{\prime }-\rho (z_{\zeta ,j})=a_{\zeta ,j}e\mbox{.}
$$
Then applying $\rho $ to the equation (\ref{zzeta}) and subtracting the
result from equation (\ref{zprime}), we obtain 
$$
\begin{array}{c}
q_{\zeta ,m}(z_{\zeta ,m+r+1}^{\prime }-\rho (z_{\zeta ,m+r+1}))=(z_{\zeta
,m+r}^{\prime }-\rho (z_{\zeta ,m+r}^{}))+ \\ 
\sum_{\ell <r}d_{\zeta ,m}^\ell (z_{\zeta ,\ell }^{\prime }-\rho (z_{\zeta
,\ell }^{}))+\sum_{k=1}^n(\psi (g_{\zeta ,m}^k)-\rho (g_{\zeta
,m}^k))+c_\zeta (m)e 
\end{array}
$$

\noindent from which, comparing coefficients in ${\Bbb Z}e$, we get
$$
q_{\zeta ,m}a_{\zeta ,m+r+1}=a_{\zeta ,m+r}+\sum_{\ell <r}d_{\zeta ,m}^\ell
a_{\zeta ,\ell }+\sum_{k=1}^nf(\varphi _\zeta ^k(m))+c_\zeta (m)\mbox{.}
$$
$\Box $



\section{(B) implies (B+)}

Now we are going to move from one combinatorial property, (B), to a stronger
one, (B+), which will allow us to construct Whitehead groups that are
strongly $\lambda $-free. Recall that in section 2 we defined the
reshuffling property (Definition \ref{reshuff}).

\begin{theorem}
\label{BtoB+}Suppose that for some regular uncountable cardinal $\lambda $,
there is a Whitehead $\lambda $-system. Then the following also holds:

(B+) there exist integers $n>0$ and $r\geq 0$, and:

\begin{enumerate}
\item  a $\lambda $-system $\Lambda =(S,\lambda _\eta ,B_\eta :\eta \in S)$
of height $n$ ;

\item  one-one functions $\varphi _\zeta ^k$ $(\zeta \in S_f$, $1\leq k\leq
n)$ with $\limfunc{dom}(\varphi _\zeta ^k)=\omega $;

\item  primes $q_{\zeta ,m}$ ($\zeta \in S_f$, $m\in \omega $);

\item  integers $d_{\zeta ,m}^\ell $ ($\zeta \in S_f$, $m\in \omega $, $\ell
<r$)
\end{enumerate}

\noindent
satisfying (a) and (b) as in (B), with the additional properties:

\begin{itemize}
\item  ${\cal S}=\{s_\eta :\eta \in S_f\}$ has the reshuffling property; and

\item  for all $\zeta \in S_f$ and $k,i \in \omega $, $\limfunc{rge}(\varphi
_\zeta ^i)\cap \limfunc{rge}(\varphi _\zeta ^k)=\emptyset $ if $i\neq k
$.
\end{itemize}
\end{theorem}

{\sc proof. }We shall refer to the data in (B+), with the given properties,
as a {\it strong Whitehead $\lambda $-system.}

Suppose that $\Lambda ^{\prime }=(S^{\prime },\lambda _\eta ^{\prime
},B_\eta ^{\prime }:\eta \in S^{\prime })$, $\varphi _\zeta ^{\prime
k},q_{\zeta ,m}^{\prime }$, and $d_{\zeta ,m}^{\prime \ell }$ is a Whitehead 
$\lambda $-system (as in (B)); in particular, ${\cal S^{\prime }}=\{s_\zeta
^{\prime }:\zeta \in S_f^{\prime }\}$ is a family of countable sets based on 
$\Lambda ^{\prime }$ , where $s_\zeta ^{\prime }=\bigcup_{k=1}^n\limfunc{rge}%
(\varphi _\zeta ^{\prime k})$. In \cite[\S VII.3A]{EM} is contained a proof
that if there exists a family ${\cal S^{\prime }}$ of countable sets based
on a $\lambda $-system $\Lambda ^{\prime }$ which is $\lambda $-free, then
there is a family ${\cal S}$ of countable sets based on a $\lambda $-system $%
\Lambda $ which has the reshuffling property. Our task is to examine the
proof and show how the transformations carried out in the proof can be done
in such a way that the additional data and properties given in (B) ---
namely the existence of the functions $\varphi _\zeta ^k$, primes $q_{\zeta
,m}$, and integers $d_{\zeta ,m}^\ell $ satisfying (b) --- continue to
hold. The transformations in question change the given ${\cal S^{\prime }}$
and $\Lambda ^{\prime }$ into ${\cal S}$ and $\Lambda $ which are {\it %
beautiful, }that is, they satisfy the following six properties:

\begin{enumerate}
\item  for $\eta $, $\nu \in S$, if $B_\eta \cap B_\nu \neq \emptyset $,
then there are $\tau \in S$ and $\alpha $, $\beta $ so that $\eta =\tau
\frown \langle \alpha \rangle $ and $\nu =\tau \frown \langle \beta \rangle $%
;

\item  for $\zeta $, $\nu \in S_f$ and $k$, $i\in \omega $, if $s_\zeta
^k\cap s_\nu ^i\neq \emptyset $ then $k=i$, $\ell (\zeta )=n=\ell (\nu )$
for some $n$ and for all $j\neq k-1$,%
\footnote{Note that this corrects an
error in \cite{EM}. A list of errata for \cite{EM} is available from the first author.}%
$\ \ \zeta (j)=\nu (j)$;

\item  for each $k$ and $\zeta$, $s_\zeta ^k$ is infinite and has a tree
structure; that is, for each $\zeta $ there is an enumeration $t_0^{k\zeta
},t_1^{k\zeta },\ldots $ of $s_\zeta ^k$ so that for all $\nu $, $\zeta \in
S_f$ and $n\in \omega $, if $t_{n+1}^{k\zeta }\in s_\nu ^k$, then $%
t_n^{k\zeta }\in s_\nu ^k$;

\item  ${\cal {S}}$ is $\lambda $-free;

\item  for all $\alpha \in E_\emptyset $, $\Lambda ^{\langle \alpha \rangle }
$ and ${\cal {S}}^{\langle \alpha \rangle }$ are beautiful; and

\item  one of the following three possibilities holds: 

\smallskip

 (a) every $\gamma \in E_\emptyset $ has cofinality $\omega $ and
there is an increasing sequence of ordinals $\{\gamma _n\colon n \in
\omega \}$ approaching $\gamma $ such that for all $\zeta \in S_f{ }$
if $\zeta (0) = \gamma $ then $s^1_\zeta = \{\langle \gamma _n,
t_n\rangle \colon n \in \omega \}$ for some $t_n$'s;  moreover,
these enumerations of the $s^1_\zeta $ satisfy the tree property of
(iii);

\smallskip
 
(b) there is an uncountable cardinal $\kappa $ and an integer  $m >
0$  so that for all  $\gamma  \in  E_\emptyset $ the cofinality of
$\gamma $ is $\kappa $ and for all  $\zeta  \in  S_f$, $\lambda
_{\zeta \rest  m} = \kappa $; moreover, for each  $\gamma  \in
E_\emptyset $ there is a strictly increasing continuous sequence
$\{\gamma _\rho \colon  \rho  < \kappa \}$ cofinal in $\gamma $
such that for all $\zeta  \in  S_f$  if  $\zeta (0) = \gamma $
then  $s^1_\zeta = \{\gamma _{\zeta (m)}\} \times  X_\zeta$  for some
$X_\zeta$;

\smallskip
 
(c) each  $\gamma  \in  E_\emptyset $ is a regular cardinal and  
$\lambda _{\langle \gamma \rangle } = \gamma $; moreover, for every 
$\zeta \in S_f$,  $s^1_\zeta  =
\{\zeta (1)\} \times  X_\zeta$ for some $X_\zeta.$

\end{enumerate}

By \cite[Thm. VII.3A.6]{EM}, if ${\cal S}$ and $\Lambda $ are beautiful,
then ${\cal S}$ has the reshuffling property. Thus it is enough to show that
we can transform ${\cal S^{\prime }}$ and $\Lambda ^{\prime }$ into a
beautiful ${\cal S}$ and $\Lambda $ and at the same time preserve the
additional structure of (B).

Let us begin with property (i). We do not change $S^{\prime }$ (the tree), but
for every $\tau \in S^{\prime } \setminus S_f$ and $\alpha \in E_\tau ^{\prime }$, we
replace $B_{\tau ,\alpha }^{\prime }$ with $B_{\tau ,\alpha }^{\prime
}\times \{\tau \}$. Define 
$$
\varphi _\zeta ^k(m)=\langle \varphi _\zeta ^{\prime k}(m),\zeta
\rest k-1\rangle \in B_{\zeta \rest k}^{\prime }\times \{\zeta \rest k-1\}.
$$
The definitions of the rest of the data are unchanged. Then (B)(b)
continues to hold since given the $c_\zeta $, define $f$ by $f^{}(\varphi
_\zeta ^k(m))=f^{\prime }(\varphi _\zeta ^{\prime k}(m))$, where $f^{\prime
} $ is the function associated with the original data (and the same $c_\zeta 
$). The function $f$ is well-defined because if $\varphi _{\zeta
_1}^k(m_1)=\varphi _{\zeta _2}^k(m_2)$, then $\varphi _{\zeta _1}^{\prime
k}(m_1)=\varphi _{\zeta _2}^{\prime k}(m_2)$. Note that property (i) implies
that $\limfunc{rge}(\varphi _\zeta ^i)\cap \limfunc{rge}(\varphi _\zeta
^k)=\emptyset $ if $i\neq k$.

Property (ii) of the definition of beautiful is handled similarly.

To obtain property (iii), we do not change $S^{\prime }$, but for all $\tau \in
S^{\prime }$ we replace $B_\tau ^{\prime }$ with $^{<\omega }B_\tau ^{\prime
}$, the set of all finite sequences of elements of $B_\tau ^{\prime }$.
Enumerate $s_\zeta ^{\prime k}$ as $\{x_{\zeta ,j}^k:j\in \omega \}$. If $%
\varphi _\zeta ^{\prime k}(m)=x_{\zeta ,j_m}^k$, define%
$$
\varphi _\zeta ^k(m)=\langle x_{\zeta ,i}^k:i\leq j_m\rangle \in
{}^{<\omega }B_{\zeta \rest k}^{\prime }\mbox{.}
$$
Given $c_\zeta $ ($\zeta \in S_f$), define $f(\varphi _\zeta
^k(m))=f^{\prime }(\varphi _\zeta ^{\prime k}(m))$, where $f^{\prime }$ is
the function associated with the original data (and the same $c_\zeta $).
Again, $f$ is well-defined.

So we can suppose that $\Lambda ^{\prime }=(S^{\prime },\lambda _\eta
^{\prime },B_\eta ^{\prime }:\eta \in S^{\prime })$ and ${\cal S^{\prime }}%
=\{s_\zeta ^{\prime }:\zeta \in S_f^{\prime }\}$ satisfy also properties 
(i),
(ii) and (iii). The proof of \cite[Thm. VII.3A.5]{EM} shows that one can define $%
\Lambda $ and ${\cal S\ }$ which are beautiful and such that

\begin{quote}
there is a one-one order-preserving map $\psi $ of $S$ into $S^{\prime }$
such that for all $\eta \in S$, $\lambda $$_\eta =\lambda^{\prime}_{\psi (\eta )}$;
and for each $\zeta \in S_f$, there is a level-preserving bijection $\theta
_\zeta :s_\zeta \rightarrow s_{\psi (\zeta )}^{\prime }$ such that for all $%
\zeta ,\nu \in S_f$, if $x\in s_{\psi (\zeta )}^{\prime }$, $y\in s_{\psi
(\nu )}^{\prime }$ and $x\neq y$, then $\theta _\zeta ^{-1}(x)\neq \theta
_\nu ^{-1}(y)$.%
\footnote{Note that this is a clarification and
correction of the first paragraph of the proof of \cite[Thm.
VII.3A.5, p. 213]{EM}. Also, in the third paragraph of that 
proof, $\psi $ should be $\psi
^{-1}$.}
\end{quote}

Observe that $\eta \in $$S_f$ if and only if $\psi (\eta )\in S_f^{\prime }$
since $\lambda _\eta =\lambda^{\prime} _{\psi (\eta )}$. We use the functions $\psi $
and $\theta _\zeta $ to define the additional data in (B): let $q_{\zeta
,m}=q_{\psi (\zeta ),m}^{\prime }$ and $d_{\zeta ,m}^\ell =d_{\psi (\zeta
),m}^{\prime \ell }$; moreover, define $\varphi _\zeta ^k(m)=$ $\theta
_\zeta ^{-1}(\varphi _{\psi (\zeta )}^{\prime k}(m))$. Given $c_\zeta $ for $%
\zeta \in S_f$, define $c_{\psi (\zeta )}^{\prime }=c_\zeta $ and let $c_\nu 
$ be arbitrary for $\nu \in S_f^{\prime }\setminus \psi [S]$. Then since the
original data satisfy (B), $f^{\prime }:\cup {\cal S^{\prime }}\rightarrow 
{\Bbb Z}$ and $a_{\nu ,j}^{\prime }$ ($\nu \in S_f^{\prime }$, $j\in \omega $%
) exist. Let $f(\varphi _\zeta ^k(m))=f^{\prime }(\varphi _{\psi (\zeta
)}^{\prime k}(m))$; the (contrapositive of the) final hypothesis on $\theta
_\zeta $ implies that $f$ is well-defined. Let $a_{\zeta ,m}=a_{\psi (\zeta
),m}^{\prime }$. Then for each $\zeta \in S_f$, the equation 
$$
q_{\psi (\zeta ),m}^{\prime }a_{\psi (\zeta ),m+r+1}^{\prime }=a_{\psi
(\zeta ),m+r}^{\prime }+\sum_{\ell <r}d_{\psi (\zeta ),m}^{\prime \ell
}a_{\psi (\zeta ),\ell }^{\prime }+\sum_{k=1}^nf^{\prime }(\varphi _{\psi
(\zeta )}^{\prime k}(m))+c_{\psi (\zeta )}^{\prime }(m) 
$$
{\it is} the desired equation 
$$
q_{\zeta ,m}a_{\zeta ,m+r+1}=a_{\zeta ,m+r}+\sum_{\ell <r}d_{\zeta ,m}^\ell
a_{\zeta ,\ell }+\sum_{k=1}^nf(\varphi _\zeta ^k(m))+c_\zeta (m)\mbox{.}
$$
$\Box $



\section{(B+) IMPLIES (A+)}

\begin{theorem}
Let $\lambda $ be a regular uncountable cardinal such that (B+) holds, i.e.,
there is a strong Whitehead $\lambda $-system. Then

(A+) there are $2^\lambda $ strongly $\lambda $-free Whitehead groups of
cardinality $\lambda $.
\end{theorem}

{\sc proof. }Given a strong Whitehead $\lambda $-system $(S,\lambda _\eta
,B_\eta :\eta \in S)$ together with $\varphi _\zeta ^k,q_{\zeta ,m},d_{\zeta ,m}^\ell $%
, we use them to define a group $G$ in terms of generators and relations.
Our group $G$ will be the group $F/K$ where $F$ is the free abelian group
with basis 
\begin{align}
\label{basis}\bigcup {\cal S}\cup \{z_{\zeta ,j}:\zeta \in S_f\mbox{, }j\in
\omega \} 
\end{align}
and $K$ is the subgroup of $F$ generated by the elements $w_{\zeta ,m}=$

\begin{align}
\label{whm}q_{\zeta ,m}z_{\zeta ,m+r+1}-z_{\zeta ,m+r}-\sum_{\ell
<r}d_{\zeta ,m}^\ell z_{\zeta ,\ell }-\sum_{k=1}^n\varphi _\zeta ^k(m)
\end{align}
for all $m\in \omega $, and $\zeta \in S_f$. Let us see first that $G$ is a
Whitehead group. (For this we need only (B).) It suffices to show that every
group homomorphism $\psi :K\longrightarrow {\Bbb Z}$ extends to a
homomorphism from $F$ to ${\Bbb Z}$. (See, for example, \cite[p.8]{EM}.)
Given $\psi $, define $c_\zeta (m)=\psi (w_{\zeta ,m})$ for all $m\in \omega 
$, and $\zeta \in S_f$. Then by (B)(b), there are integers $a_{\zeta ,j}$ ($%
\zeta \in S_f,$ $j\in \omega $) and a function $f:\bigcup {\cal 
S}\to 
 {\Bbb Z }$ such that for all $\zeta \in S_f$ and $m\in \omega $,
\begin{align} \label{qa}c_\zeta (m)=q_{\zeta ,m}a_{\zeta
,m+r+1}-a_{\zeta ,m+r}-\sum_{\ell \ <r}d_{\zeta ,m}^\ell a_{\zeta ,\ell
}-\sum_{k=1}^nf(\varphi _\zeta ^k(m)). \end{align} Define $\theta
:F\longrightarrow {\Bbb Z}$ by setting $\theta \rest \bigcup {\cal S}
=f$ and $\theta (z_{\zeta ,j})=a_{\zeta ,j}$. We just
need to check that $\theta $ extends $\psi .$ But for all $\zeta \in
S_f$ and $ m\in \omega $, we have $$ \theta (w_{\zeta
,m})=q_{\zeta,m}a_{\zeta ,m+r+1}-a_{\zeta ,m+r}-\sum_{\ell <r}d_{\zeta
,m}^\ell a_{\zeta ,\ell }-\sum_{k=1}^nf(\varphi _\zeta ^k(m)) $$ by the
definitions of $\theta $ and of $w_{\zeta ,m}$. Thus $\theta (w_{\zeta
,m})=c_\zeta (m)=\psi (w_{\zeta ,m})$, by (\ref{qa}).

Next let us show that $G$ is not free. (Here again, we need only (B).) The
proof is essentially the same as that of Lemma VII.3.9 of \cite[pp. 205f]{EM}%
, but we will give a somewhat different version of the proof here. The proof
proceeds by induction on $n$ where $n$ is the height of our $\lambda $%
-system. For each $\alpha <\lambda $, let $G_\alpha $ be the subgroup of $G$
generated by 
$$
\{z_{\zeta ,j}:\zeta \in S_f,\zeta (0)<\alpha \mbox{, }j\in \omega \}\cup
\bigcup \{s_\zeta :\zeta \in S_f,\zeta (0)<\alpha \}. 
$$
It suffices to prove that for all $\alpha $ in a stationary subset of $%
\lambda ,$ $G_{\alpha +1}/G_\alpha $ is not free (cf. \cite[IV.1.7]{EM}). In
fact, we will show that $G_{\alpha +1}/G_\alpha $ is not free when $\alpha $
is a limit point of $E_\emptyset $ and belongs to $C\cap E_\emptyset $,
where $C$ is the cub

\smallskip

\begin{quote}
$\{\alpha <\lambda :$ whenever $\varphi _\zeta ^1(m)\in \bigcup
\{B_{\left\langle \beta \right\rangle }:\beta <\alpha \}$ then $\exists
\sigma \in S_{f\mbox{ }}^{}$ with $\sigma (0)<\alpha $ and $\varphi _\zeta
^1(m)\in \limfunc{rge}(\varphi _\sigma ^1)\}.$
\end{quote}
\smallskip

We begin with the case $n=1$. Then for all $\alpha \in C\cap E_\emptyset $
such that $\alpha $ is a limit point of $E_\emptyset $, $G_{\alpha +1}/G_a$
is non-free because it is as described in the first part of Lemma \ref
{nonfree} (with generators $\{z_{\left\langle \alpha \right\rangle ,j}:j\in
\omega \}$), since for all $m\in \omega $, $\varphi _{\left\langle \alpha
\right\rangle }^1(m)\in B_{\left\langle \alpha \right\rangle }=\bigcup
\{B_{\left\langle \beta \right\rangle }:\beta <\alpha \}$ by the definition
of a $\lambda $-system (because $\alpha $ is a limit point of $E_\emptyset $%
) and hence $\varphi _{\left\langle \alpha \right\rangle }^1(m)\in G_\alpha $
since $\alpha \in C$.

Now suppose $n>1$ and the result is proved for $n-1$. Again, let $\alpha \in
C\cap E_\emptyset $ such that $\alpha $ is a limit point of $E_\emptyset $.
Again we have that $\varphi _\zeta ^1(m)\in G_\alpha $ for all $m\in \omega $
when $\zeta (0)=\alpha $. We will consider the $\lambda _{\left\langle
\alpha \right\rangle }$-system $\Lambda ^{\left\langle \alpha \right\rangle
} $. (See Definition \ref{genlsys}.) Note that $\Lambda ^{\left\langle
\alpha \right\rangle }$ has height $n-1$, and the group $G_{\alpha
+1}/G_\alpha $ is defined as in (\ref{basis}) and (\ref{whm}) relative to
this $\lambda _{\left\langle \alpha \right\rangle }$-system. Hence by
induction $G_{\alpha +1}/G_\alpha $ is not free.

Finally, we will use the reshuffling property given by (B+) to prove that $G$
is strongly $\lambda $-free. As in the proof of \cite[VII.3.11]{EM}, we will
prove that for all $\alpha \in \lambda \cup \{-1\}$ and all $\beta $$>\alpha 
$, $G_\beta /G_{\alpha +1}$ is free. Let $I=\{\zeta \in S_f\colon \zeta
(0)<\beta \}$, and let $<_I$ be the well-ordering given by the reshuffling
property for $I$ and $\alpha $. Let $s_\zeta ^k$ denote $\limfunc{rge}%
(\varphi _\zeta ^k)$. We claim that there is a basis ${\cal Z}_{\beta
,\alpha }$ of $G_\beta /G_{\alpha +1}$ consisting of the cosets of the
members of the following two sets:%
$$
\begin{array}{l}
\{z_{\zeta ,j}:\alpha <\zeta (0)<\beta , 
\mbox{ and either }j<r\mbox{ or}\  \\ \exists k\mbox{ s.t. }\varphi _\zeta
^k(j-r)\notin \bigcup \{s_\nu ^k:\nu <_I\zeta \}\} 
\end{array}
$$

\noindent and 
$$
\begin{array}{l}
\{\varphi _\zeta ^k(m):\varphi _\zeta ^k(m)\notin \bigcup \{s_\nu ^k:\nu
<_I\zeta \} 
\mbox{ and} \\ \exists i<k[\varphi _\zeta ^i(m)\notin \bigcup \{s_\nu ^i:\nu
<_I\zeta \}]\}\mbox{.}
\end{array}
$$

To see that the elements of ${\cal Z}_{\beta ,\alpha }$ generate $G_\beta
/G_{\alpha +1}$, we proceed by induction with respect to $<_I$ to show that
the coset of every $z_{\zeta ,j}$ ($\zeta (0)<\beta ,j\in \omega $) and the
coset of every element of $s_\zeta $ ($\zeta (0)<\beta $) is a linear
combination of the elements of ${\cal Z}_{\beta ,\alpha }$. Since $s_\zeta
\setminus \bigcup \{s_\nu ^{}:\nu <_I\zeta \}$ is infinite, for each $j\in
\omega $ such that $z_{\zeta ,j}+G_{a+1}\notin {\cal Z}_{\beta ,\alpha }$,
there is $t>j$ such that $z_{\zeta ,t}+G_{\alpha +1}$ belongs to ${\cal Z}%
_{\beta ,\alpha }$. Without loss of generality, $t=j+1.$ Then 
$$
z_{\zeta ,j}=q_{\zeta ,j-r}z_{\zeta ,t}-\sum_{\ell <r}d_{\zeta ,j-r}^\ell
z_{\zeta ,\ell }-\sum_{k=1}^n\varphi _\zeta ^k(j-r) 
$$
by (\ref{whm}). By induction each $\varphi _\zeta ^k(j-r)+G_{a+1}$ is a
linear combination of members of ${\cal Z}_{\beta ,\alpha }$ (because $%
\varphi _\zeta ^k(j-r)\in \bigcup \{s_\nu ^k:\nu <_I\zeta \}$ since $%
z_{\zeta ,j}+G_{a+1}\notin {\cal Z}_{\beta ,\alpha })$; hence $z_{\zeta
,j}+G_{a+1}$ belongs to the subgroup generated by the members of ${\cal Z}%
_{\beta ,\alpha }$.

For each $m,i\in \omega $, if $\varphi _\zeta ^i(m)\in \bigcup \{s_\nu
^k:\nu <_I\zeta \}$, then by induction $\varphi _\zeta ^i(m)+G_{\alpha +1}$
is a linear combination of elements of ${\cal Z}_{\beta ,\alpha }$.
Otherwise, $\varphi _\zeta ^i(m)+G_{\alpha +1}$ belongs to ${\cal Z}_{\beta
,\alpha }$ unless $i$ is minimal such that $\varphi _\zeta ^i(m)\notin
\bigcup \{s_\nu ^k:\nu <_I\zeta \}$. But in the latter case,%
$$
\sum_{k=1}^n\varphi _\zeta ^k(m)=q_{\zeta ,m}z_{\zeta ,m+r+1}-z_{\zeta
,m+r}-\sum_{\ell <r}d_{\zeta ,m}^\ell z_{\zeta ,\ell } 
$$
so its coset is a linear combination of elements of ${\cal Z}_{\beta ,\alpha
}$; thus since $\varphi _\zeta ^k(m)+G_{\alpha +1}\in {\cal Z}_{\beta
,\alpha }$ for all $k\neq i$, $\varphi _\zeta ^i(m)+G_{\alpha +1}$ belongs
to the subgroup generated by the elements of ${\cal Z}_{\beta ,\alpha }$.
This completes the proof that ${\cal Z}_{\beta ,\alpha }$ is a generating
set. To see that the elements of ${\cal Z}_{\beta ,\alpha }$ are
independent, compare coefficients in $F$.

To construct not just one but 2$^\lambda $ different strongly $\lambda $%
-free Whitehead groups, we use a standard trick: write $E_\emptyset $ as the
disjoint union $\coprod_{\sigma <\lambda }X_\sigma $ of $\lambda $
stationary sets; then for every non-empty subset $W$ of $\lambda $, do the
construction above for the generalized $\lambda $-system $\Lambda
=(S_W,\lambda _\zeta ,B_\zeta :\zeta \in S_W)$ with $E_\emptyset =$ $%
\coprod_{\sigma \in W}X_\sigma $, i.e., where $S_W=\{\zeta \in S:\zeta
(0)\in \coprod_{\sigma \in W}X_\sigma \}$. $\Box $



\section{Appendix: Non-free Whitehead  implies 2-uniformization}

A {\it ladder system} on a stationary subset $E$ of $\omega _1$ is an
indexed family of functions $\{\eta $$_\alpha :\alpha \in E\}$ such that
each $\eta _\alpha :\omega \rightarrow \alpha $ is strictly increasing and $%
\sup (\limfunc{rge}(\eta _\alpha ))=\alpha $. If $\{\varphi _\alpha :\alpha
\in I\}$ is an indexed family of functions each with domain $\omega $, we
say that it has the $2${\it -uniformization property} provided that for
every family of functions $c_\alpha :\omega \rightarrow 2=\{0,1\}$ ($\alpha
\in I$), there exists a function $H$ such that for all $\alpha \in I$, $%
H(\varphi _\alpha (n))$ is defined and equals $c_\alpha (n)$ for all but
finitely many $n$. It is not hard, given a ladder system on $E$ which has
the $2$-uniformization property, to construct, explicitly (by generators and
relations), a non-free Whitehead group. (See \cite[Prop. XII.3.6]{EM}.) It
is more difficult to go the other way: starting with an arbitrary non-free
Whitehead group of cardinality $\aleph_1$ to show that there exists a ladder system on a stationary
subset of $\omega _1$ which has the $2$-uniformization property. This was
left to the reader in the original paper by the second author 
\cite[Thm. 3.9, p. 277]{Sh80}. The only published proof is a rather
complicated one in \cite[\S XII.3]{EM}; so considering the importance of
this result, it seems to us worthwhile to give another proof which is
conceptually and technically simpler than that one. The proof given here
resembles the original proof found by the second author, which was also the
basis of the proofs in \cite{EMS} and in this paper.

Our goal is to prove the following.

\begin{theorem}
\label{appmain}If there is a non-free Whitehead group $A$ \ of cardinality $%
\aleph _1$, \ then there is a ladder system $\{\eta _\alpha :\alpha \in E\}$
on a stationary set $E$ which has the 2-uniformization property.
\end{theorem}

We begin with an observation. It is sufficient to show that the hypothesis
of Theorem \ref{appmain} implies that there is a family $\{\varphi _\alpha
:\alpha \in E\}$ of functions which has the 2-uniformization property and is 
{\it based on an $\omega _1$-filtration}, that is, indexed by a stationary
subset $E$ of $\omega _1$ and such that there is a continuous ascending
chain $\{B_\nu :\nu \in \omega \}$ of countable sets such that for all $%
\alpha \in E$, $\varphi _\alpha :\omega \rightarrow B_\alpha $. (Note that
what we are talking about, in the language of the preceding sections, is a
family of countable sets based on an $\aleph _1$-system.) Indeed, by a
suitable coding we can assume that $B_\alpha =\alpha $ and if the range of $%
\varphi _\alpha $ is not cofinal in $\alpha $, we can choose a ladder $\eta
_\alpha ^{\prime }$ on $\alpha $, replace $\varphi _\alpha (n)$ by $%
\left\langle \varphi _\alpha (n),\eta _\alpha ^{\prime }(n)\right\rangle $,
and re-code, to obtain a ladder system on $E\cap C$, (for some cub $C$)
which has the $2$-uniformization property.

\relax Fromnow on, let $A$ denote a non-free Whitehead group of cardinality $%
\aleph _1$. Then we can write $A$ as the union, $A=\cup _{\nu <\omega
_1}A_\nu $, of a continuous chain of countable free 
subgroups; since $A$ is not
free, we can assume that there is a stationary subset $E$ of $\omega _1$%
(consisting of limit ordinals) such that for all $\alpha \in E$ $A_{\alpha
+1}/A_\alpha $ is not free. By Pontryagin's Criterion we can assume without
loss of generality that $A_{\alpha +1}/A_\alpha $ is of finite rank and, in
fact, that every subgroup of $A_{\alpha +1}/A_\alpha $ of smaller rank is
free. Since

\begin{quote}
($\star $){\it \ whenever $E=\cup _{n\in \omega }E_n$, at least one of the $%
E_n$ is stationary}
\end{quote}

\noindent (cf. \cite[Cor. II.4.5]{EM}) we can also assume that all of the $%
A_{\alpha +1}/A_\alpha $ (for $\alpha \in E$) have the same rank $r+1$ ($%
r\geq 0$). In order to make clear the ideas involved in the proof of the
Theorem, we will give the proof first in the special case when $r=0$, i.e., $%
A_{\alpha +1}/A_\alpha $ is a rank one non-free group when $\alpha \in E$,
and then describe how to handle the extension to the general case. In fact
this special case divides into two subcases: using ($\star $) and replacing $%
A_{\alpha +1}$ by a subgroup if necessary, we can assume that either

\begin{enumerate}
\item  for all $\alpha \in E$, $A_{\alpha +1}/A_\alpha $ has a type all of
whose entries are $0$'s or $1$'s [and there are infinitely many $1$'s]; or

\item  there is a prime $p$ such that for all $\alpha \in E$, the type of $%
A_{\alpha +1}/A_\alpha $ is \\ $(0,0,...0,\infty ,0,...)$ where the $\infty $
occurs in the $p$th place.
\end{enumerate}

\noindent (See \cite[pp. 107ff]{F}.) We next give the easy combinatorial
lemmas needed for the first, and simplest, subcase.

\begin{lemma}
\label{1}Suppose $Y$ and $Y^{\prime }$ \ are finite subsets of an abelian
group $G$ such that $|Y|^2<|Y^{\prime }|$. Then there exists $b\in Y^{\prime
}$ such that $Y$ \ and $b+Y$ are disjoint. [Here $b+Y=\{b+y:y\in Y\}$.]
\end{lemma}

{\sc proof. }Choose $b\in Y^{\prime }\setminus \{x-y:x,y\in Y\}$. $\Box $

\begin{lemma}
\label{2}For any positive integer $p>1$ there are integers $a_0$ and $a_1$
and a function $F_p:{\Bbb Z}/p{\Bbb Z}\rightarrow 2=\{0,1\}$ such that for all
$m\in {\Bbb Z}$ with $(2|m|+1)^2<p$, $F_p(m+a_\ell +p{\Bbb Z})=\ell $,
for $\ell =0,1$.
\end{lemma}

{\sc proof.} Let $a_0=0$ and let $a_1=b$ as in Lemma \ref{1}, where $G={\Bbb Z%
}/p{\Bbb Z}=Y^{\prime }$ and $Y=\{m+p{\Bbb Z}:$ $(2|m|+1)^2<p\}$. Then since $%
Y=a_0+Y$ and $a_1+Y$ are disjoint, we can define $F_p$. (Note that $F_p$ is
a set function, not a homomorphism.) $\Box $

 \  

{\sc proof of theorem \ref{appmain} }({\it in special subcase} (i)): For all $%
\alpha \in E$ there is an infinite set $P_\alpha $ of primes such that 
$$
A_{\alpha +1}/A_\alpha \cong \{\frac mn\in {\Bbb Q}:n\mbox{ is a product of
{\it distinct} primes from}\ P_\alpha \}\mbox{.}
$$
Then if $P_\alpha =\{p_{\alpha ,n}:n\in \omega \}$, $A_{\alpha +1}$ is
generated over $A_\alpha $ by a subset $\{y_{\alpha ,n}:n\in \omega \}$
satisfying the relations (and only the relations) 
$$
(\dagger )\mbox{ }\vspace*{.3in}\mbox{ }p_{\alpha ,n}y_{\alpha
,n+1}=y_{\alpha ,0}^{}-g_{\alpha ,n} 
$$
for some $g_{\alpha ,n}\in A_\alpha $. We define $\varphi _\alpha
(n)=\left\langle p_{\alpha ,n},g_{\alpha ,n}\right\rangle $. Then $\{\varphi
_\alpha :\alpha \in E\}$ is based on an $\omega _1$-filtration, in fact on
the chain $\{{\Bbb Z}\times A_\alpha :\alpha <\omega _1\}$.

Given functions $c_\alpha :\omega $ $\rightarrow 2$, we are going to define
a homomorphism $\pi :A^{\prime }\rightarrow A$ with kernel ${\Bbb Z}e$ and
then use the splitting $\rho :A\rightarrow A^{\prime }$ to define the
uniformizing function $H$.

We define $\pi _\nu :A_\nu ^{\prime }\rightarrow A_\nu $ inductively along
with a set function $\psi _\nu $$:A_\nu \rightarrow A_\nu ^{\prime }$ such
that $\pi _\nu \circ \psi _\nu =1_{A_\nu }$. The crucial case is when $\pi
_\alpha $ and $\psi _\alpha $ have been defined and $\alpha \in E$. (When $%
\alpha \notin E$ we can use the fact that $\limfunc{Ext}(A_{\alpha +1},
{\Bbb Z)}\rightarrow \limfunc{Ext}(A_\alpha ,{\Bbb Z)}$ is onto.) We define $%
A_{\alpha +1}^{\prime }$ by generators $\{y_{\alpha ,n}^{\prime ^{}}:n\in
\omega \}$ over $A_\alpha ^{\prime }$ satisfying relations 
$$
(\dagger \dagger )\mbox{ }\vspace*{.3in}\mbox{ }p_{\alpha ,n}y_{\alpha
,n+1}^{\prime }=y_{\alpha ,0}^{\prime }-\psi _\alpha (g_{\alpha ,n})+a_\ell
e 
$$
where $a_\ell $ is as in Lemma \ref{2} for $p=p_{\alpha ,n}$ and $\ell
=c_\alpha (n)$.

In the end we let $\pi =\cup _\nu \pi _\nu :A^{\prime }=\cup _\nu A_\nu
^{\prime }\rightarrow A$ and $\psi =\cup _\nu \psi _\nu $. Then since $A$ is
a Whitehead group, there exists a homomorphism $\rho $ such that $\pi \circ
\rho =1_A$. For any $g\in A$, $\psi (g)-\rho (g)\in \ker (\pi )={\Bbb Z}e$;
we will abuse notation and identify $\psi (g)-\rho (g)$ with the unique
integer $k$ such that $\psi (g)-\rho (g)=ke$. For any $w\in \cup _{\alpha
\in E}\limfunc{rge}(\varphi _\alpha )$, if $w=\left\langle p,g\right\rangle $%
, let $H(w)=F_p(\psi (g)-\rho (g)+p{\Bbb Z})$.

Note that $w$ may equal $\varphi _\alpha (n)$ ($=\left\langle p_{\alpha
,n},g_{\alpha ,n}\right\rangle $) for many pairs $(\alpha ,n)$. To see that
this definition of $H$ works, fix $\alpha \in E$. For any $n\in \omega $,
applying $\rho $ to equation ($\dagger $) and subtracting from equation ($%
\dagger \dagger $), we have%
$$
p_{\alpha ,n}(y_{\alpha ,n+1}^{\prime }-\rho (y_{\alpha ,n+1}))=y_{\alpha
,0}^{\prime }-\rho (y_{\alpha ,0})-(\psi _\alpha (g_{\alpha ,n})-\rho
(g_{\alpha ,n}))+a_\ell 
$$
so that $\psi _\alpha (g_{\alpha ,n})-\rho (g_{\alpha ,n})$ is congruent to $%
y_{\alpha ,0}^{\prime }-\rho (y_{\alpha ,0})+a_\ell $ mod $p_{\alpha ,n}$.
Then if $n$ is large enough, $(2|$ $y_{\alpha ,0}^{\prime }-\rho (y_{\alpha
,0})|+1)^2<p_{\alpha ,n}$ so by choice of $F_{p_{\alpha ,n}}$ and $a_\ell $, 
$$
\begin{array}{c}
H(\varphi _\alpha (n))=F_{p_{\alpha ,n}}(\psi _\alpha (g_{\alpha ,n})-\rho
(g_{\alpha ,n})+p_{\alpha ,n} 
{\Bbb Z})= \\ F_{p_{\alpha ,n}}(y_{\alpha ,0}^{\prime }-\rho (y_{\alpha
,0})+a_\ell +p_{\alpha ,n}{\Bbb Z})=\ell =c_\alpha (n)\mbox{.}
\end{array}
$$
This completes the proof in the first special subcase.

For the purposes of the second special subcase we need another combinatorial
lemma.

\begin{lemma}
\label{sub2}Fix a positive integer $p>1$. Define a strictly increasing
sequence of positive integers $t_i$ inductively, as follows. Let $t_0=0$. If 
$t_{i-1}$ has been defined for some $i\geq 1$, let $t_i=t_{i-1}+d_i$ where $%
d_i$ is the least positive integer such that $%
(2p^{t_{i-1}}+1)^2p^{2t_{i-1}}<p^{d_i}.$ Then for every $i\geq 1$ there
exists a function 
$$
F_i:{\Bbb Z}/p^{t_i}{\Bbb Z}\rightarrow 2
$$
and integers $a_n^\ell \in \{0,...,p-1\}$ ($t_{i-1}\leq n<t_i,\ell =0,1$)
such that whenever $|m_0|\leq p^{t_{i-1}}$ and $a_j\in \{0,...,p-1\}$ for $%
j<t_{i-1}$, then for $\ell =0,1$%
$$
F_i(m_0+\sum_{j<t_{i-1}}p^ja_j+\sum_{n=t_{i-1}}^{t_i-1}p^na_n^\ell +p^{t_i}%
{\Bbb Z})=\ell \mbox{.}
$$
\end{lemma}

{\sc proof.} We apply Lemma \ref{1} to the sets $Y=\{m_0+%
\sum_{j<t_{i-1}}p^ja_j+p^{t_i}Z:|m_0|\leq p^{t_{i-1}}$, $a_j\in
\{0,...,p-1\}\}$ (which has cardinality $\leq (2p^{t_{i-1}}+1)p^{t_{i-1}}$)
and $Y^{\prime }=\{\sum_{n=t_{i-1}}^{t_i-1}p^nx_n+p^{t_i}Z:x_n\in
\{0,...,p-1\}\}$ (which has cardinality $p^{d_i}$), to get $b\in Y^{\prime }$%
. Then choose $a_n^0=0$ for all $n$, and $a_n^1$ so that $%
\sum_{n=t_{i-1}}^{t_i-1}p^na_n^1=b$ and define $F_i$ as in Lemma \ref{2}. $%
\Box $

\  

{\sc proof of theorem \ref{appmain} }({\it in special subcase} (ii)): For all $%
\alpha \in E$%
$$
A_{\alpha +1}/A_\alpha \cong \{\frac mn\in {\Bbb Q}:n\mbox{ is a power of }%
p\} 
$$
Then $A_{\alpha +1}$ is generated over $A_\alpha $ by a subset $\{y_{\alpha
,n}:n\in \omega \}$ satisfying the relations (and only the relations) 
$$
(\dagger )\mbox{ }\vspace*{.3in}\mbox{ }py_{\alpha ,n+1}=y_{\alpha
,n}-g_{\alpha ,n} 
$$
for some $g_{\alpha ,n}\in A_\alpha $. Let $\varphi _\alpha (m)=\left\langle
g_{\alpha ,j}:j<t_{m+1}\right\rangle $ for all $m\in \omega .$ Given
functions $c_\alpha :\omega $ $\rightarrow 2$, we define $\pi _\nu :A_\nu
^{\prime }\rightarrow A_\nu $ with kernel ${\Bbb Z}e$ inductively along with
a set function $\psi _\nu $$:A_\nu \rightarrow A_\nu ^{\prime }$ such that $%
\pi _\nu \circ \psi _\nu =1_{A_\nu }$. The crucial case is when $\pi _\alpha 
$ and $\psi _\alpha $ have been defined and $\alpha \in E$. Then we define $%
A_{\alpha +1}^{\prime }$ by generators $\{y_{\alpha ,n}^{\prime ^{}}:n\in
\omega \}$ over $A_\alpha ^{\prime }$ satisfying relations 
$$
(\dagger \dagger )\mbox{ }\vspace*{.3in}\mbox{ }py_{\alpha ,n+1}^{\prime
}=y_{\alpha ,n}^{\prime }-\psi _\alpha (g_{\alpha ,n})+a_n^{\ell (n)}e 
$$
where $\ell (n)$ is taken to be $c_\alpha (i-1)$ when $t_{i-1}\leq n<t_i$.
In the end we let $\pi =\cup _\nu \pi _\nu :A^{\prime }=\cup _\nu A_\nu
^{\prime }\rightarrow A$ and $\psi =\cup _\nu \psi _\nu $. Then since $A$ is
a Whitehead group, there exists a homomorphism $\rho $ such that $\pi \circ
\rho =1_A$. For any $w\in \cup _{\alpha \in E}\limfunc{rge}(\varphi _\alpha
) $, if $w=\left\langle g_j:j<t_i\right\rangle $, let $H(w)=F_i(%
\sum_{n<t_i}^{}p^n(\psi (g_n)-\rho (g_n))+p^{t_i}{\Bbb Z}{)}$. To see that
this works, fix $\alpha \in E$ and for $i\geq 1$, consider $w_i=\varphi
_\alpha (i-1)=\left\langle g_{\alpha ,j}:j<t_i\right\rangle $. From the
equations ($\dagger $), for $n\leq t_i$ we obtain that 
$$
p^{t_i}y_{\alpha ,t_i}=y_{\alpha ,0}-\sum_{n<t_i}^{}p^ng_{\alpha ,n} 
$$
If we apply $\rho $ to this and subtract from the corresponding equation
derived from ($\dagger \dagger $) we obtain that $\sum_{n<t_i}^{}p^n(\psi
(g_{\alpha ,n})-\rho (g_{\alpha ,n}))$ is congruent to 
$$
(y_{\alpha ,0}^{\prime }-\rho (y_{\alpha ,0}))+\sum_{n<t_i}^{}p^na_n^{\ell
(n)}\ 
$$
mod $p^{t_i}$. So if $|y_{\alpha ,0}^{\prime }-\rho (y_{\alpha ,0})|\leq
p^{t_{i-1}\mbox{ }}$, then by our choice of $F_i$ and the $a_n^{\ell (n)}$
for $t_{i-1}\leq n<t_i$, $H(w_i)$ equals $c_\alpha (i-1)$.

 This completes
the proof of Theorem \ref{appmain} when $r=0$.

\  

{\sc proof of theorem \ref{appmain} }({\it in the general case}): In the
general case without loss of generality we have either

\begin{enumerate}
\item  for all $\alpha \in E$, $A_{\alpha +1}/A_\alpha $ has a free subgroup 
$L_\alpha /A_\alpha $ of rank $r$ such that $A_{\alpha +1}/L_\alpha $ has a
type all of whose entries are $0$'s or $1$'s [and there are infinitely many $%
1$'s]; or

\item  there is a prime $p$ such that for all $\alpha \in E$, $A_{\alpha
+1}/A_\alpha $ has a free subgroup $L_\alpha /A_\alpha $ of rank $r$ such
that the type of  $A_{\alpha +1}/L_\alpha $ is $(0,0,...0,\infty ,0,...)$ where the 
$\infty $ occurs in the $p$th place.
\end{enumerate}

In other words, $A_{\alpha +1}$ is generated by $A_\alpha $ and a subset $%
\{z_{\alpha ,k}:k=1,...,r\}\cup $ $\{y_{\alpha ,n}:n\in \omega \}$ modulo
(only) the relations in $A_\alpha $ plus relations:

\begin{enumerate}
\item  ($\dagger $) $p_{\alpha ,n}y_{\alpha ,n+1}=y_{\alpha
,0}+\sum_{k=1}^r\mu _{\alpha ,k}(n)z_{\alpha ,k}-g_{\alpha ,n}$ for some
family of distinct primes $p_{\alpha ,n}$ and $\mu _{\alpha ,k}(n)\in {\Bbb Z}
$ , $g_{\alpha ,n}\in A_\alpha $; or

\item  ($\dagger $) $py_{\alpha ,n+1}=y_{\alpha ,n}+\sum_{k=1}^r\mu _{\alpha
,k}(n)z_{\alpha ,k}-g_{\alpha ,n}$ for some $\mu _{\alpha ,k}(n)\in {\Bbb Z}$
and $g_{\alpha ,n}\in A_\alpha $ for each $n\in \omega $.
\end{enumerate}

\ 

For use in (the harder) subcase (ii), define a strictly increasing sequence of
positive integers $t_i$ inductively, as follows. Let $t_0=0$. If $t_{i-1}$
has been defined for some $i\geq 1$, let $t_i=t_{i-1}+d_i$ where $d_i$ is
the least positive integer such that%
$$
(2p^{t_{i-1}}+1)^{2r+2}p^{2t_{i-1}}<p^{d_i}\mbox{.}
$$
Then we have the following generalization of Lemma \ref{sub2}. (Note that
when $r=0$ the sequence $\mu $ is empty.)

\begin{lemma}
\label{3}Fix $p>1$ and $r\geq 0$. For every sequence of functions $\mu
=\left\langle \mu _1,...,\mu _r\right\rangle $, where $\mu _k:\omega
\rightarrow {\Bbb Z}$ and every $i\geq 1$ there exists a function
$$
F_{i,\mu }:{\Bbb Z}/p^{t_i}{\Bbb Z}\rightarrow 2
$$
and integers $a_{n,\mu }^\ell \in \{0,...,p-1\}$ ($t_{i-1}\leq n<t_i,\ell
=0,1$) such that $F_{i,\mu }$ and $a_{n,\mu }^\ell $ depend only on $\mu
\rest t_i $ ($=\left\langle \mu _1\rest t_i,...,\mu _r\rest t_i\right\rangle $) and are
such that whenever $m_0,...,m_r$ are integers with $|m_k|\leq p^{t_{i-1}}$
for all $k\leq r$ and $a_j\in \{0,...,p-1\}$ for $j<t_{i-1}$, then%
$$
F_{i,\mu }(m_0+\sum_{k=1}^r(\sum_{j<t_i}^{}p^j\mu
_k(j))m_k+\sum_{j<t_{i-1}}p^ja_j+\sum_{n=t_{i-1}}^{t_i-1}p^na_{n,\mu }^\ell
+p^{t_i}{\Bbb Z})=\ell \mbox{.}
$$
\end{lemma}

{\sc proof. } We apply Lemma \ref{1} with $G={\Bbb Z}/p^{t_i}{\Bbb Z}$,
$$
\begin{array}{c}
Y=\{m_0+\sum_{k=1}^r(\sum_{j<t_i}^{}p^j\mu
_k(j))m_k+\sum_{j<t_{i-1}}p^ja_j+p^{t_i} 
{\Bbb Z}: \\ |m_k|\leq p^{t_i}\mbox{ , for all }k\leq r\mbox{, }a_j\in
\{0,...,p-1\}\} 
\end{array}
$$
and%
$$
Y^{\prime }=\{\sum_{n=t_{i-1}}^{t_i-1}p^nx_n+p^{t_i}{\Bbb Z}:x_n\in
\{0,...,p-1\}\}\mbox{.}
$$
and proceed as in the proof of Lemma \ref{sub2}. $\Box $

Similarly we have the following generalization of Lemma \ref{2} for use in
subcase (i).

\begin{lemma}
\label{4}Given $p>1$ and $r\geq 0$, and a sequence of integers $\mu
=\left\langle \mu _1,...,\mu _r\right\rangle $, let $t_p$ be maximal such
that $(2t_p+1)^{2r+2}<p$. Then there exists a function 
$$
F_{p,\mu }:{\Bbb Z}/p{\Bbb Z}\rightarrow 2
$$
and integers $a_{p,\mu }^\ell \in \{0,...,p-1\}$ ($\ell =0,1$) such that
whenever $m_0,...,m_r$ are integers such that $|m_k|\leq t_p$ for all $k\leq
r$, then%
$$
F_{p,\mu }(m_0+\sum_{k=1}^r\mu _km_k+a_{p,\mu }^\ell +p{\Bbb Z})=\ell \mbox{.}
$$
$\Box $
\end{lemma}

Now define the function $\varphi _\alpha $ with domain $\omega $ by letting

\begin{enumerate}
\item  $\varphi _\alpha (m)=\left\langle \left\langle \mu _{\alpha
,k}(m):k=1,...,r\right\rangle ,p_{\alpha ,m},g_{\alpha ,m}\right\rangle $; or

\item  $\varphi _\alpha (m)=\left\langle \left\langle \mu _{\alpha
,k}(n):k=1,...,r\right\rangle ,g_{\alpha ,n}:n<t_{m+1}\right\rangle $.
\end{enumerate}

Given functions $c_\alpha :\omega ${\it $\rightarrow 2$}, we define $\pi
_\nu :A_\nu ^{\prime }\rightarrow A_\nu $ inductively along with a set
function $\psi _\nu $$:A_\nu \rightarrow A_\nu ^{\prime }$ such that $\pi
_\nu \circ \psi _\nu =1_{A_\nu }$. The crucial case is when $\pi _\alpha $
and $\psi _\alpha $ have been defined and $\alpha \in E$. Then we define $%
A_{\alpha +1}^{\prime }$ by generators $\{z_{\alpha ,k}^{\prime
}:k=1,...,r\}\cup \{y_{\alpha ,n}^{\prime ^{}}:n\in \omega \}$ over $%
A_\alpha ^{\prime }$ satisfying relations

\begin{enumerate}
\item  ($\dagger \dagger $) $p_{\alpha ,n}y_{\alpha ,n+1}^{\prime
}=y_{\alpha ,0}^{\prime }+\sum_{k=1}^r\mu _{\alpha ,k}(n)z_{\alpha
,k}^{\prime }-\psi _\alpha (g_{\alpha ,n})+a_{p_{\alpha ,n},\mu }^\ell e$; or

\item  ($\dagger \dagger $) $py_{\alpha ,n+1}^{\prime }=y_{\alpha
,n}^{\prime }+\sum_{k=1}^r\mu _{\alpha ,k}(n)z_{\alpha ,k}^{\prime }-\psi
_\alpha (g_{\alpha ,n})+a_{\alpha ,n,\mu }^\ell e$
\end{enumerate}

\noindent where $a_{p_{\alpha ,n},\mu }^\ell $ (respectively, $a_{\alpha
,n,\mu }^\ell $) is as in Lemma \ref{4} (respectively, Lemma \ref{3}) for $%
\ell =c_\alpha (n)$ (respectively, $\ell =c_\alpha (i-1)$ if $t_{i-1}\leq
n<t_i$ $)$ (and the appropriate prime or primes are used).

In the end we use a splitting $\rho $ of $\pi =\cup _\nu \pi _\nu :A^{\prime
}=\cup _\nu A_\nu ^{\prime }\rightarrow A$ to define $H(w)$ as follows:

\begin{enumerate}
\item  if $w=\left\langle \left\langle \mu _k:k=1,...,r\right\rangle
,p,g\right\rangle $, let $H(w)=F_{p,\mu }(\psi (g)-\rho (g)+p{\Bbb Z})$; or

\item  if $w=\left\langle \left\langle \mu _k(n):k=1,...,r\right\rangle
,g_n:n<t_i\right\rangle $, let $H(w)= \newline F_{i,\mu }(\sum_{n<t_i}^{}p^n(\psi
(g_n)-\rho (g_n))+p^{t_i}{\Bbb Z})$.
\end{enumerate}

\noindent Then we check as before that this definition works. $\Box $




\begin{thebibliography}{99}
\bibitem{E}  P. C. Eklof, {\it Infinitary equivalence of abelian groups},
Fund. Math. {\bf 81} (1974), 305--314.

\bibitem{EM}  P. C. Eklof and A. H. Mekler, {\bf {Almost Free Modules}},
North-Holland (1990).

\bibitem{EMS}  P. C. Eklof, A. H. Mekler and S. Shelah, {\it Uniformization
and the diversity of Whitehead groups}, Israel J. Math {\bf 80} (1992),
301-321.

\bibitem{ES}  P. C. Eklof and S. Shelah, {\it On Whitehead Modules}, J.
Algebra {\bf 142} (1991), 492--510.

\bibitem{F}  L. Fuchs, {\bf Infinite Abelian Groups}, vol II, Academic Press
(1973).

\bibitem{MS}  M. Magidor and S. Shelah, {\it When does almost free imply
free? (For groups, transversal etc.)}, to appear in Journal Amer. Math. Soc.

\bibitem{Sh74}  S. Shelah, {\it Infinite abelian groups, Whitehead problem
and some constructions}, Israel J. Math {\bf {18}} (1974), 243--25.

\bibitem{Sh75}  S. Shelah, {\it A compactness theorem for singular
cardinals, free algebras, Whitehead problem and transversals}, Israel J.
Math., {\bf 21}, 319--349.

\bibitem{Sh80}  S.{\ Shelah, }{\it Whitehead groups may not be free even
assuming CH, II}, Israel J. Math. {\bf 35} (1980), 257--285.

\bibitem{Sh85}  S. Shelah, {\it {Incompactness in regular cardinals}}, Notre
Dame J. Formal Logic {\bf {26 }}(1985), 195--228.

\bibitem{T}  J. Trlifaj, {\it Non-perfect rings and a theorem of Eklof and
Shelah,} Comment. Math. Univ. Carolinae {\bf 32} (1991), 27--32.
\end{thebibliography}
\end{document}